\documentclass[final,onefignum,onetabnum]{siamart171218}
\usepackage[utf8]{inputenc}
\usepackage{amsfonts}
\usepackage{amsopn} 
\usepackage{mathtools} 
\usepackage{bm} 
\usepackage{commath} 
\usepackage{stmaryrd} 
\usepackage{nicefrac}
\usepackage{textcomp}

\usepackage{multirow} 
\usepackage{colortbl} 
\newcommand{\myrowcolor}{\rowcolor[gray]{0.925}}
\usepackage{booktabs} 
\usepackage[font=footnotesize,position=b]{subcaption}
\usepackage{adjustbox}

\usepackage{tikz,tikzscale} 
\usetikzlibrary{calc}
\usetikzlibrary{shadings}
\usetikzlibrary{patterns}
\definecolor{apply}{rgb}{0.3,.7,0.} 
\definecolor{applyface}{rgb}{.9,.95,0.}
\definecolor{colorA}{rgb}{1.0,0.6,0.28} 
\definecolor{colorB}{rgb}{0.,0.,0.} 

\definecolor{slategrey}{rgb}{0.44,0.50,0.56}
\definecolor{darkorange}{rgb}{0.93,0.6,0.}
\usetikzlibrary{patterns}
\usepackage{graphicx}
\graphicspath{{grids/},{reductions/},{throughput/}}

\usepackage{algpseudocode}

\usepackage[color=gray,backgroundcolor=gray,
            ]{todonotes}
\reversemarginpar

\newsiamremark{remark}{Remark}
\newsiamremark{hypothesis}{Hypothesis}
\Crefname{hypothesis}{Hypothesis}{Hypotheses}
\newsiamthm{claim}{Claim}
\newsiamremark{assumption}{Assumption}

\let\domain\Omega
\newcommand{\mapping}[1]{F_{#1}}
\newcommand{\tria}[1]{\mathcal{T}_{#1}}

\newcommand{\cell}{K}
\newcommand{\refcell}{\hat{\cell}}
\newcommand{\face}{e}
\newcommand{\poly}[1]{\mathbb P_{#1}}
\newcommand{\tppoly}[1]{\mathbb{Q}_{#1}}

\newcommand{\avg}[1]{\{ #1 \}}
\newcommand{\normal}{\bm{n}}

\newcommand{\dgedges}[1]{\mathcal{E}_{#1}}
\newcommand{\dgedgesi}[1]{\mathcal{E}_{#1}^{\circ}}
\newcommand{\dgedgesb}[1]{\mathcal{E}_{#1}^{\partial}}
\newcommand{\ippre}{\hat\gamma}
\newcommand{\ip}{\gamma}
\newcommand{\interval}{I}
\newcommand{\relaxpre}{\omega}

\newcommand{\qweight}{\omega}
\newcommand{\tensoru}{U}
\newcommand{\tensorJxW}{T}
\newcommand{\multiindex}[1]{{#1_1,\dots,#1_d}} 
\newcommand{\mulindex}[1]{{#1_1,\dots,#1_d}} 
\newcommand{\grad}{\nabla}
\newcommand{\gradref}{\hat\nabla}
\newcommand{\doftransfer}{\Pi_{\ell,\cell}}

\DeclareMathOperator{\diag}{diag}

\newcommand{\xref}{\hat x}
\newcommand{\bfx}{\bm x}
\newcommand{\bfxref}{\bm{\xref}}
\newcommand{\dx}{\,dx}

\newcommand{\dbfx}{\,d\bfx}

\newcommand{\dsbfx}{\,d\sigma(\bfx)}

\newcommand{\bigo}{\mathcal{O}}

\newcommand{\jacobianr}{\nabla F}
\DeclareMathOperator{\determinant}{det}

\newcommand{\hadap}{\circ}
\newcommand{\resred}{\delta_{\textit{red}}}
\newcommand{\restol}{\varepsilon_{\textit{tol}}}
\newcommand{\fracsteps}{\nu_{\textit{frac}}}
\newcommand{\normalizationcmplx}{n_{\textit{sub}}}
\newcommand{\ccmplx}{C_\textit{cmplx}}
\newcommand{\ordercmplx}{\textit{order}}

\headers{Fast Tensor Product Schwarz Smoothers}{J. Witte, D. Arndt, and G. Kanschat}

\title{Fast Tensor Product Schwarz Smoothers for High-Order Discontinuous Galerkin Methods\thanks{
\funding{The research presented in this article was funded by the German Research Foundation (DFG) under
the project “High-order discontinuous Galerkin for the EXA-scale” (ExaDG) within the priority program
“Software for Exascale Computing” (SPPEXA).\newline 
This manuscript has been authored by UT-Battelle, LLC under Contract No.
  DE-AC05-00OR22725 with the U.S. Department of Energy. The United States Government
  retains and the publisher, by accepting the article for publication, acknowledges
  that the United States Government retains a non-exclusive, paid-up, irrevocable,
  worldwide license to publish or reproduce the published form of this manuscript,
  or allow others to do so, for United States Government purposes. The Department
  of Energy will provide public access to these results of federally sponsored
  research in accordance with the DOE Public Access Plan
  (http://energy.gov/downloads/doe-public-access-plan).
}}}

\author{Julius Witte\thanks{Interdisziplin\"ares Zentrum f\"ur Wissenschaftliches Rechnen (IWR), Universität Heidelberg, Im Neuenheimer Feld 205, 69120 Heidelberg. Germany.\newline Email: \email{julius.witte/guido.kanschat@iwr.uni-heidelberg.de}}
\and Daniel Arndt\thanks{Computational Engineering and Energy Sciences Group,
         Oak Ridge National Laboratory;
         Oak Ridge, TN, USA.\newline Email: \email{arndtd@ornl.gov}}
\and Guido Kanschat\footnotemark[2]}

\begin{document}

\maketitle

\begin{abstract}
  In this article, we discuss the efficient implementation of powerful domain decomposition smoothers for multigrid methods for high order discontinuous Galerkin (DG) finite element methods. In particular, we study the inversion of matrices associated to mesh cells and to the patches around a vertex, respectively, in order to obtain fast local solvers for additive and multiplicative subspace correction methods. The effort of inverting local matrices for tensor product polynomials of degree $k$ is reduced from $\mathcal O(k^{3d})$ to  $\mathcal O(dk^{d+1})$ by exploiting the separability of the differential operator and resulting low rank representation of its inverse as a prototype for more general low rank representations.
\end{abstract}



\section{Introduction}

This article shows that powerful multigrid smoothers based on domain decomposition with cells and vertex patches as subdomains can be implemented very efficiently using fast diagonalization.
In particular, we show that now, as matrix-free application of operators associated with finite element bilinear forms is state of the art, implementation of powerful smoothers can be accomplished with the same asymptotic complexity with respect to polynomial degree.
The technique demonstrated for the Laplacian can be applied to any separable operator.

Multigrid and domain decomposition methods are the two classes of solvers or preconditioners which allow the solution of discretized elliptic partial differential equations with linear complexity in the number of degrees of freedom. This is true at least for mesh refinement.
When the polynomial degree is increased, point smoothers on the fine mesh must be replaced by more complex methods with in general superlinear complexity.
For instance, nonoverlapping subdomain smoothers using the inversion of cell matrices of the interior penalty method (IP) have yielded very effective multigrid methods for higher order discretizations of the Laplacian~\cite{Kanschat08smoother}, reaction-diffusion systems~\cite{LuceroKanschat18}, or radiation transport~\cite{Lucero18}.
For divergence constrained problems, they are not sufficient and we have to resort to overlapping patches of $2^d$ hypercube cells around a single vertex.
Then, we obtain effective multigrid methods for divergence-dominated problems~\cite{ArnoldFalkWinther97Hdiv}, the Stokes problem~\cite{KanschatMao15}, or a Darcy-Stokes-Brinkman system~\cite{KanschatLazarovMao17}.
While these methods are very effective in the sense of few iteration steps, computation time can become unfeasible in a standard implementation, if large cell matrices are inverted with an algorithm of cubic complexity. 

Modern hardware favors algorithms performing complex operations on small data sets, since memory access is by far more expensive in terms of time and energy than computation.
Thus, it was observed for instance in~\cite{KronbichlerKormann12,kronbichler2019fast,MuethingPiatkowskiBastian17} that implementations based on stored sparse matrices, which have a computational intensity of one FLOP for each entry read from memory are not competitive.
On the other hand, once the computational intensity is high enough that computation dominates memory access, it is not only worthwhile, but mandatory to optimize the computational part of algorithms.
This has been achieved for application of finite element operators, where most codes now prefer integration of bilinear forms over mesh cells computed on the fly to stored matrices.
While unfeasible with a standard quadrature with complexity $k^6$ in the polynomial degree $k$ in three dimensions, applications of a local matrix to a local vector can be performed at low arithmetic cost and complexity of order $k^4$ using the technique of sum factorization. This technique has been introduced in the context of spectral methods in~\cite{orszag80}.

For effective multigrid smoothers, we need the solution of local problems in addition to operator application.
Hence, we turn to low rank tensor representations of the local matrices and their inverses to yield a similar reduction of complexity, following the idea in~\cite{lynch64}.
There, a low rank tensor representation of the inverse of a separable (discrete) operator on a tensor product mesh is presented.
We apply it as a local solver on overlapping subdomains on Cartesian meshes and as an approximate local solver on more general geometries.

Kronecker decompositions of separable operators have been used as preconditioners in~\cite{Beuchler03,lottes05,stiller2017,stiller2016}. In~\cite{Beuchler03}, the one dimensional local problems are preconditioned by a wavelet method and then used in a block preconditioner of the global system, splitting the bubble degrees of freedom on edges and the interior of cells, respectively, from those those in vertices for continuous finite elements.
Methods closer to ours were introduced in the context of continuous finite elements in~\cite{lottes05,stiller2017}. There, a cell based smoother is introduced which augments each cell by a few layers of support points with associated basis functions from neighboring cells. In~\cite{stiller2016}, a similar smoother for discontinuous Galerkin methods is presented and compared to an augmented patch of two cells sharing a face.
From the point of view of data structures, these approaches are more complicated than ours, since we use cell-wise data instead of augmenting by neighboring shape functions.
The block preconditioners in~\cite{pazner18} based on Kronecker decompositions are not restricted to separable operators. However, the successive Kronecker singular value decompositions~\cite{vanloan93, vanloan00} of the local solvers on cells require $\bigo(k^{2d-1})$ instead of $\bigo(k^{d+1})$ arithmetic operations per element, losing optimality in three dimensions.
An alternative approach to efficient local solvers based on localized matrix-free methods is the iterative solution of the local problems up to a fixed accuracy, see~\cite{BastianMuellerMuethingPiatkowski18}.

Algorithms in this article have been developed with vectorizing multicore architectures and their high cost for memory access compared to computation in mind. Here, we study their computational complexity only, showing that even with purely sequential arithmetic, we obtain good multigrid performance. The reason are low iteration counts combined with an implementation of the smoothers with low overhead.
Parallel implementation, where the balance between computation and memory access becomes important, is deferred to a forthcoming study.

This article is organized as follows: in the following section, we introduce the model problem, its discretization by the interior penalty method, and the multilevel Schwarz methods we use for preconditioning together with some results on convergence speed. In \Cref{sec:tpelements}, we present the efficient implementation of these smoothers for separable operators and results for their computational effort. \Cref{sec:unstructured} discusses the application of inexact local solvers on more general meshes and their impact on performance. Our findings are summarized in \Cref{sec:conclusions}.


\section{Multilevel interior penalty methods}
\label{sec:MLIP}
In this work, we discuss a method for the model problem of Poisson's equation
\begin{equation}
    \label{eq:model}
    \begin{aligned}
    - \Delta u &= f &&\text{ in } \Omega,\\
    u &= g &&\text{ on } \partial\Omega,
    \end{aligned}
\end{equation}
where $\Omega$ is a polygonal domain in $\mathbb{R}^d$ with $d = 2, 3$.
$f$ and $g$ are given functions in $L^2(\Omega)$ and $L^2(\partial\Omega)$, respectively. We point out that we used the Laplacian as a simple example, but that it can be replaced by any separable operator. In case of nonsymmetric operators, eigenvalues below must be replaced by singular values.

\subsection{Discretization}
\label{sec:disc}
The model problem is discretized by means of the symmetric interior penalty method (SIPG) following \cite{ArnoldBrezziCockburnMarini02,Arnold82}.
To this end, we subdivide the domain $\Omega$ into meshes $\tria \ell$ for levels $\ell=0,\dots,L$, where the finest level $L$ is the actual discretization level on which we want to solve and the intermediate levels $\ell<L$ form the hierarchy for the geometric multigrid method.
Each mesh consists of a collection of quadrilateral/hexahedral cells $\cell$, which are obtained by a mapping $\mapping{\cell}$ from the reference cell $\refcell = [0,1]^d$.
The relation of these meshes is defined by induction as follows: starting from a coarse mesh $\tria{0}$ consisting of few cells at most, we generate a hierarchy of meshes $\tria{\ell}$ for levels $\ell = 0,\ldots,L$ by recursively splitting each cell in $\tria{\ell}$ with respect to its midpoint into $2^d$ children  in $\tria{\ell+1}$. These meshes are nested in the fashion that every cell of $\tria{\ell}$ is equal to the union of its $2^d$ children in $\tria{l+1}$ as well as conforming in the sense that either any edge/face of the cell is at the domain's boundary or a complete edge/face of an adjacent cell.

The shape function space $V(\refcell)$ on the reference cell consists of the discontinuous, tensor product polynomials $\tppoly{k}$. Its basis $\{\hat\varphi_i\}$ consists of tensor products of Lagrangian interpolation polynomials of degree $k$ with respect to the Gauss-Lobatto points on the reference interval. More details will be provided in \Cref{sec:tpelements}.
The shape function spaces $V(\cell)$ on an actual grid cell are obtained by composition with the mapping $F_\cell$ such that $\varphi_{\cell,i}(\bfx) = \hat\varphi_i(F_{\cell}^{-1} \bfx)$. The finite element space on level $\ell$ is then defined by
\begin{equation}
    \label{eq:discFEspace}
    V_\ell := \{ v \in L^2(\domain) \mid v|_{\cell} \in V(\cell) \: \forall \cell \in \tria \ell \} = \bigoplus_{\forall \cell \in \tria \ell} V(\cell).
\end{equation}
The indexing of the basis $\{\varphi_{\cell,i}\}$ follows the structure as a direct sum.
This basis defines by duality the coefficient space $\mathbb R^{n_\ell}$ of the same dimension as $V_\ell$ equipped with the Euclidean inner product. In computations, this is the inner product used to compute norms, such that we will identify $V_\ell$ with the coefficient space and do not distinguish them in notation.

Let $\dgedgesi \ell$ be the set of all interior interfaces between two cells $\cell^{\pm}$. Then, we refer to traces of functions $v\in V_\ell$ on $\face\in\dgedgesi\ell$ taken from cell $\cell^\pm$ as $v^\pm$, respectively. For such a function, we define the ``averaging operator''
\begin{equation}
    \label{eq:jumpavgscalar}
    \avg{v}(\bfx) = \frac{1}{\sqrt2} \bigl(v^+(\bfx) + v^-(\bfx)\bigr),
    \qquad \bfx\in \face.
\end{equation}
On a face at the boundary, denoted by $\face\in\dgedgesb\ell$, there is only a trace from the interior and thus we define
\begin{equation}
    \label{eq:avgbdry}
    \avg{v}(\bfx) = v(\bfx),
    \qquad \bfx\in \face.
\end{equation}
Using $\normal$ as the outward pointing normal of the cell $\cell$ at face $\face$, we introduce the interior penalty bilinear form
\begin{align}
\begin{aligned}
    \label{eq:bfsipg}
    a_\ell (u, v) :=
    & \int_{\tria\ell} \nabla u \cdot \nabla v \dbfx \\
    & + \int_{\dgedges\ell} \left( 
    \ip_\face \, \avg{u \normal} \cdot \avg{v \normal}
    - \avg{\nabla u} \cdot \avg{v \normal} 
    - \avg{u \normal} \cdot \avg{\nabla v}
    \right) \dsbfx.
\end{aligned}  
\end{align}
Here, the integrals over sets of cells (faces) are to be understood as the sum of the individual integrals over cells (faces).
From left to right we refer to the four integrals on the right hand side of~\eqref{eq:bfsipg} as the \emph{bulk}, \emph{penalty}, \emph{consistency} and \emph{adjoint consistency} term.
We still have to define the edge-wise penalty parameter $\ip_\face$, which penalizes the jumps of the solution and yields stability of the bilinear form \cite[\S2.2.8]{Kanschat03habil}. It is of the form
\begin{equation}
  \label{eq:penalty}
  \ip_\face = \ippre\; k(k+1) \left(\frac{1}{h^+} + \frac{1}{h^-}\right) \quad \text{on } \dgedgesi{h},
\end{equation}
where $h^\pm$ is the (average) length of cell $K^\pm$ orthogonal to the common edge $e$. On boundary edges $e \in \dgedgesb{h}$ we let $h^+ = h^- = h$, where $h$ is the length of the corresponding cell orthogonal to $e$.
The factor $\ippre$ is chosen equal to one on Cartesian elements and has to be increased on non-Cartesian elements to preserve stability of the discretization.
Finally, we can state the interior penalty discretization of the model problem~\eqref{eq:model}:
find $u_L \in V_L$ such that
\begin{equation}
    a_L (u_L, v) = 
    \int_{\Omega} f v \dbfx 
    + \int_{\partial\Omega} \left( \ip_\face \, g v - g \normal \cdot \nabla v \right) \dsbfx 
    \qquad \forall v \in V_L.
    \label{eq:sipg}
\end{equation}

\subsection{Geometric Multigrid}
\label{sec:multigrid}

We follow~\cite{GopalakrishnanKanschat03} in the definition of the
geometric multigrid algorithm for the interior penalty method. More
precisely, we state the V-cycle algorithm used for preconditioning in \Cref{alg:gmg}.
\begin{algorithm}[tp]
\caption{V-Cycle on level $\ell$}\label{alg:gmg}
\begin{algorithmic}[1]
  \Procedure{MG$_{\ell}$}{$x_{\ell},b_{\ell}$}
  \If{$\ell=0$}
  \State{$x_0 \gets A_0^{-1} b_0$} \Comment{coarse grid solver}
  \EndIf
\For{$k=1$ to $m_{\textit{pre}}$}
\State $x_{\ell} \gets S^{\textit{pre}}_{\ell} (x_{\ell}, b_{\ell})$
\Comment{pre-smoothing}
\EndFor
\State $b_{\ell-1} \gets I^\downarrow_{\ell-1} \bigl(b_{\ell} - A_{\ell} x_{\ell}\bigr)$ \Comment{restriction}
\State $e_{\ell-1} \gets \text{MG}_{\ell-1}(0,b_{\ell-1})$ \Comment{recursion}
\State $x_{\ell} \gets x_{\ell} +  I^\uparrow_{\ell-1} e_{\ell-1}$ \Comment{prolongation}
\For{$k=1$ to $m_{\textit{post}}$}
\State $x_{\ell} \gets S^{\textit{post}}_{\ell} (x_{\ell}, b_{\ell})$
\Comment{post-smoothing}
\EndFor
\State \Return{$x_{\ell}$} 
\EndProcedure
\end{algorithmic}
\end{algorithm}
The operators used there are as follows: $A_\ell$ refers to the level
matrix associated with the interior penalty bilinear form
in~\eqref{eq:bfsipg}. The operator $I^\uparrow_{\ell}: V_{\ell}\to V_{\ell+1}$
is the prolongation operator. Since under our assumptions the spaces
are nested, this is simply the embedding from $V_{\ell}$ into
$V_{\ell+1}$.
The restriction operator $I^\downarrow_{\ell}: V_{\ell+1}\to V_{\ell}$
is the adjoint of the prolongation operator with respect to the
Euclidean inner product in the coefficient spaces. This definition corresponds to the transpose matrix.
The operators $S^{\text{pre}}_{\ell}$ and $S^{\text{post}}_{\ell}$
are the smoothers on level $\ell$ described in detail in the
next subsection.


\subsection{Schwarz Smoothers}
\label{sec:schwarz}

We use the terms domain decomposition smoothers or Schwarz smoothers in the context of multigrid methods and many, very small subdomains, on which we solve the differential equation exactly.
Examples from the literature are the $\bm H^{\text{div}}$ and  $\bm H^{\text{curl}}$ smoothers in~\cite{ArnoldFalkWinther97Hdiv,ArnoldFalkWinther00} or cell based smoothers for the interior penalty method in~\cite{Kanschat08smoother}. The first group has been generalized successfully to Stokes~\cite{KanschatMao15} and Darcy-Stokes-Brinkman~\cite{KanschatLazarovMao17} problems. The second class has been generalized to singularly perturbed reaction-diffusion problems in~\cite{LuceroKanschat18}, where we also generalized the convergence analysis in~\cite{DryjaKrzyzanowski16} to quadrilateral and hexahedral meshes. Thus, we consider two classes of domain decomposition smoothers with local solvers on cells and vertex patches, respectively.

\begin{figure}[tp]
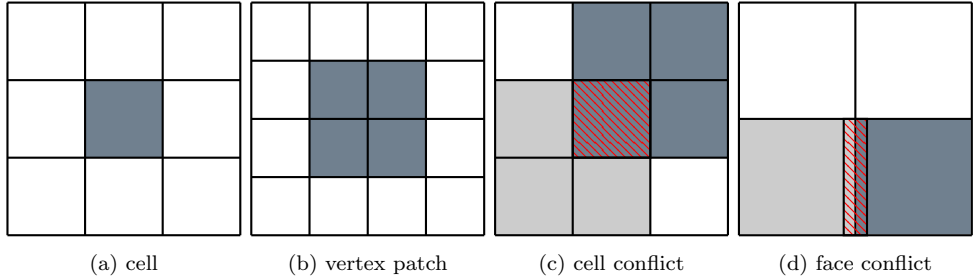
\label{fig:patches_conflicts}
\centering
\begin{subfigure}{0.24\textwidth}
	\includegraphics[width=\textwidth]{cell-dg.tikz}
	\caption{cell}
	\label{sfig:cp}
\end{subfigure}
\begin{subfigure}{0.24\textwidth}
	\includegraphics[width=\textwidth]{vertex-patch-dg.tikz}
	\caption{vertex patch}
	\label{sfig:vp}
\end{subfigure}
\begin{subfigure}{0.24\textwidth}
	\includegraphics[width=\textwidth]{vp-cell-conflict.tikz}
	\caption{cell conflict}
	\label{sfig:cell-conflict}
\end{subfigure}
\begin{subfigure}{0.24\textwidth}
	\includegraphics[width=\textwidth]{cp-face-conflict.tikz}
	\caption{face conflict}
	\label{sfig:face-conflict}
\end{subfigure}
\caption{Cell and vertex patch subdomains and conflicts due to overlap and transfer over faces}
\end{figure}

\begin{enumerate}
\item \textbf{cell based} smoothers: each subdomain of the spatial
  decomposition on level $\ell$ consists of a single cell of the mesh
  $\tria\ell$ as depicted in \Cref{sfig:cp}. After enumerating the cells in $\tria\ell$ as $\cell_j$ with $j=1,\dots,J_\ell$, the subspaces
  $V_{j;\ell} \subset V_\ell$ consist of functions with
  support in the cell $\cell_j$. As the spatial decomposition is
  nonoverlapping and we use discontinuous finite elements, $V_\ell$ is
  the disjoint union of the $V_{j;\ell}$.
\item \textbf{vertex patch} smoothers: each subdomain $\domain_j$ consists
  of all cells sharing the vertex $v_j$ of $\tria\ell$ (after
  enumeration) as shown in \Cref{sfig:vp}. The subspaces $V_{j;\ell} \subset V_\ell$ for
  $j=1,\dots,J_\ell$, where $J_\ell$ is the number of interior vertices in
  $\tria\ell$, consist of functions with support in the cell
  $\domain_j$. As typically $2^d$ cells share a vertex and a cell has
  $2^d$ vertices, the spatial decomposition is overlapping and the
  union of the subspaces is not disjoint.
\end{enumerate}

In both cases, we define the local solvers $P_{j;\ell}\colon V_\ell\to V_{j;\ell}$ by
\begin{equation}
    \label{eq:ritzproj}
    a_\ell (P_{j;\ell} u_\ell, v) = a_\ell(u_\ell, v) \quad \forall v \in V_{j;\ell}.
\end{equation}
We refer to the operator associated with the bilinear form restricted to $V_{j;\ell}$ as $A_{j;\ell}$. From now on we suppress the level index $\ell$ in expressions like $P_{j;\ell}$ and $V_{j;\ell}$.
We define the additive Schwarz smoother on level $\ell$ as
\begin{equation}
    \label{eq:additiveschwarz}
    P_{\ell;\text{ad}} := \relaxpre \sum_{j=1}^{J_\ell} P_j = \relaxpre \sum_{j=1}^{J_\ell} R^T_j A^{-1}_j R_j A_\ell,
\end{equation}
where $\relaxpre$ is a relaxation parameter. $R_j\colon V_\ell \to V_{j}$ is the restriction operator and its transpose the embedding. The form on the right highlights the structure as a product of the system matrix $A_{\ell}$ and the additive Schwarz preconditioner $A_{\ell;\textit{ad}}^{-1} = \sum_{j=1}^{J_l} R^T_j A^{-1}_j R_j$.

The multiplicative Schwarz operator, in its standard form is defined by 
\begin{equation}
  \label{eq:multschwarz1}
  P_{\ell;\text{mu}} := I - \bigl(I-\relaxpre P_{J_\ell}\bigr)\cdots
  \bigl(I- \relaxpre P_{2}\bigr) \bigl(I- \relaxpre P_{1}\bigr).
\end{equation}

While we do not evaluate parallel performance in this article, but rather focus on the numerical efficiency of the smoothing methods, namely the number of arithmetic operations needed, we nevertheless have parallel execution by vectorization, multi-threading, and MPI-parallelization on distributed systems in mind.  Neither the standard form of the multiplicative smoothers, nor the additive smoother with vertex patches are suited for such parallelism. While the additive vertex patch smoother suffers from race conditions, both the cell and the vertex patch multiplicative smoother are inherently sequential. Therefore, we use ``coloring'' of the mesh cells in order to recover potential parallelism.

Coloring refers to splitting the index set $\mathcal{J} = \{1,\dots,J\}$ of subdomains into disjoint subsets $\mathcal{J}_c$ with $c=1,\dots,C$, such that the operations within each subset can be performed in parallel without causing conflicts.

Race conditions are conflicts due to simultaneous reading and writing. They appear in the additive vertex patch smoother, if two local solvers executed in parallel are writing into the data of the same cell. Therefore, the two colored patches in \Cref{sfig:cell-conflict} may not be processed in parallel. Coloring for this algorithm is designed such that two patches of the same color do not share a common cell. For regular meshes, this can be achieved by $2^d$ parquetings of the domain with possible omission of strips at the boundary. 

\begin{figure}[tp]
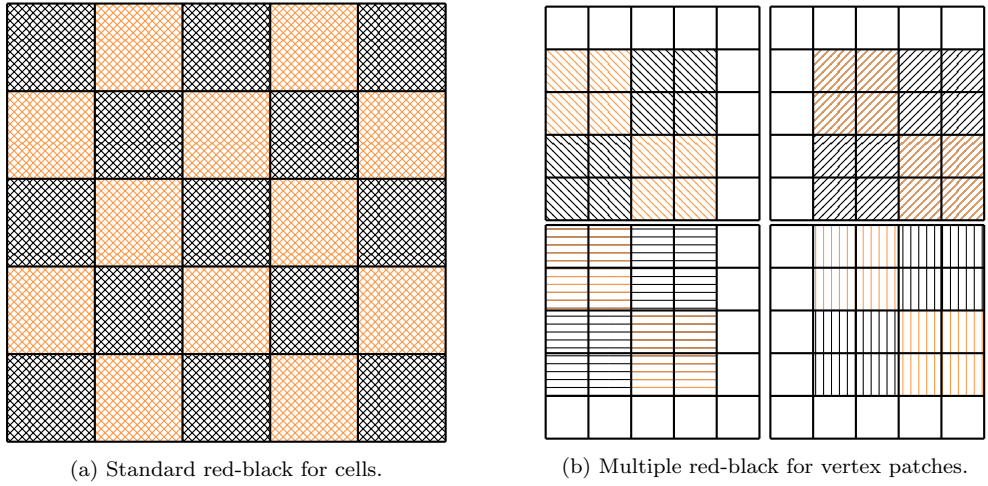

\centering
\begin{subfigure}{0.45\textwidth}
	\includegraphics[width=\textwidth]{red-black-cp.tikz}
	\caption{Standard red-black for cells.}
	\label{fig:red-black}
\end{subfigure}
\hfill
\begin{subfigure}{0.45\textwidth}
	\includegraphics[width=.49\textwidth]{coloring01-vp.tikz}
	\includegraphics[width=.49\textwidth]{coloring23-vp.tikz}

	\includegraphics[width=.49\textwidth]{coloring45-vp.tikz}
	\includegraphics[width=.49\textwidth]{coloring67-vp.tikz}
	\caption{Multiple red-black for vertex patches.}
	\label{fig:multi-red-black}
\end{subfigure}
\caption{Coloring for multiplicative algorithms.}
\label{fig:coloring}
\end{figure}

For the multiplicative algorithm, the goal of coloring is not just avoiding race conditions, it is recovering parallelism at all. To this end, we note that
\begin{align}
\begin{aligned}
    (I-P_i)(I-P_j)
    &= (I-R_i^T A_i^{-1}R_i A_\ell) (I-R_j^T A_i^{-1}R_j A_\ell)\\
    &= I - P_i - P_j + R_i^T A_i^{-1}(R_iA_\ell R_j^T) A_j^{-1}R_jA_\ell.
\end{aligned}
\end{align}
The parenthesis in the last term evaluates to zero if and only if $V_{i}$ is $A$-orthogonal to $V_{j}$. Since operator application involves face terms, $A$-orthogonality is violated if two subdomains share a common face as in \Cref{sfig:face-conflict}. It can be avoided on regular meshes by red-black coloring for the cell based smoother as in \Cref{fig:red-black}, yielding 2 colors in any dimension.

The multiplicative vertex patch smoother combines the conflicts of the additive vertex patch and the multiplicative cell based smoother. A coloring for this situation is shown in \Cref{fig:multi-red-black}, where we combine the parqueting for the vertex patch with red-black coloring into $2^{d+1}$ colors.

These simple coloring algorithms based on checkerboards and parqueting reach their limits on unstructured meshes. The finite element library \texttt{deal.II} provides graph-based coloring based on the DSATUR algorithm in~\cite{brelaz79} in a parallel version. We do not discuss its details here and refer the reader to~\cite{turcksin16}. Compared to the parqueting algorithms above, it generates more colors with smaller subsets within each color. Examples are provided below in \Cref{tab:cubemvp}.

The colored versions of the Schwarz operators read
\begin{equation}
  \label{eq:coloredschwarz}
  P_{\ell;\text{ad}} := \relaxpre \sum_{c=1}^{C} \sum_{j \in \mathcal{J}_c} P_j,
  \qquad
  P_{\ell;\text{mu}} = I - \biggl(I - \relaxpre \sum_{j \in \mathcal{J}_C} P_j\biggr)
  \cdots \biggl(I - \relaxpre \sum_{j \in \mathcal{J}_1} P_j\biggr).
\end{equation}
Each factor of the multiplicative algorithm contains sums due to $A$-orthogonality of the subspaces within the same color.
Note that the colored additive algorithm is mathematically equivalent to~\eqref{eq:additiveschwarz}, while the multiplicative version may differ due to the reordering of factors.

Since the additive Schwarz operator factors as $P_{\ell;ad} = A^{-1}_{\ell;ad} A_\ell$, the smoothing step $S_{\ell;ad}$ can be encoded as shown in \Cref{alg:addsmoother}.
\begin{algorithm}
\caption{Additive Schwarz Smoother}\label{alg:addsmoother}
\begin{algorithmic}[1]
\Procedure{$S_{\ell;\textit{ad}}$}{$x_{\ell},b_{\ell}$}
\State $r_\ell \gets b_\ell - A_\ell x_\ell$ \Comment{update residual}
\For{$c=1$ to $C$}
\State $x_\ell \gets x_\ell + \sum_{j \in \mathcal{J}_c} R^T_j A^{-1}_j R_j r_\ell $ \Comment{apply local solvers}
\EndFor  
\State \Return{$x_\ell$} 
\EndProcedure
\end{algorithmic}
\end{algorithm}
Here, we use a ``for'' loop for sequential operations, while ``$\Sigma$'' indicates a possibly parallel summation. Using the form in~\eqref{eq:coloredschwarz}, the multiplicative smoother $S_{\ell;mu}$ can be implemented in a very similar way, with the single change that the residual update happens inside the loop over all colors, as in \Cref{alg:multsmoother}.
\begin{algorithm}
\caption{Multiplicative Schwarz Smoother}\label{alg:multsmoother}
\begin{algorithmic}[1]
\Procedure{$S_{\ell;\textit{mu}}$}{$x_{\ell},b_{\ell}$}
\For{$c=1$ to $C$}
\State $r_\ell \gets b_\ell - A_\ell x_\ell$ \Comment{update residual}
\State $x_\ell \gets x_\ell + \sum_{j \in \mathcal{J}_c} R^T_j A^{-1}_j R_j r_\ell $ \Comment{apply local solvers}
\EndFor  
\State \Return{$x_\ell$} 
\EndProcedure
\end{algorithmic}
\end{algorithm}
Thus, both methods are implemented as a ``short'' product (in the sense of a sequence of operations) over all colors with parallel, additive smoothers for each color.

Since parallelization is only viable within each color, a small number of colors with many subdomains per color is desirable. This holds in particular for the multiplicative algorithm, where an operator application (residual update) is applied for each color. Therefore, whenever the meshes are regular, the optimal coloring by parqueting described above should be used. If the meshes are not regular, we fall back to the DSATUR algorithm mentioned above.

\subsection{Efficiency of the smoothers}\label{sec:robustcart}
Here, we provide numerical evidence that the smoothers discussed above yield very efficient multigrid methods in terms of iteration counts. Summarizing our findings, we show that iteration counts are independent of the mesh size $h_L$ of the finest level. Furthermore, our results indicate a slow deterioration of convergence steps for the nonoverlapping smoothers with increasing polynomial degree. The overlapping vertex patch smoother even seems to profit from higher degrees.

We present results on Cartesian meshes in two and three dimensions and refer to \Cref{sec:unstructured} for the more general case. The coarse mesh $\tria0$ is the decomposition of the square or cube $[0,1]^d$ into $2^d$ congruent cells, consequently it consists of one vertex patch. Each subsequent level is obtained by the refinement algorithm outlined in \Cref{sec:disc}. 
We use the V-cycle with a single pre- and post smoothing step as preconditioner in the conjugate gradient solver (CG) and in the generalized minimal residual method (GMRES) for the additive and multiplicative versions, respectively. The stopping criterion of the Krylov subspace methods is a relative residual reduction of $\resred = 10^{-8}$. On the coarse mesh, we solve with a relative accuracy of $10^{-8}$ using a Chebyshev solver (see~\cite{varga09}) with the additive cell based smoother as preconditioner.

Due to their efficiency, multiplicative vertex patch smoothers require three or less iterations such that we consider fractional iterations $\fracsteps$ for a more accurate assessment of their performance. If the reduction by $\resred$ is achieved after $\nu$ iterations, we compute
\begin{align}
    \fracsteps = \nu - 1 + \frac{\log(e_{\nu-1}/\restol)}{\log(e_{\nu-1}/e_\nu)},
\end{align}
where $\restol = e_0\resred$ and $e_k$ is the energy norm of the error and the Euclidean norm of the residual after $k$ steps of the CG and GMRES algorithms, respectively.
The right-hand side $f$ of our model problem~\eqref{eq:model} is manufactured such that the exact solution $u$ is given by a superposition of ``normalized'' multivariate Gaussian bell curves,
\begin{equation}
    u (\boldsymbol x) = \frac{1}{\sqrt{2\pi}\sigma} \sum_{i=1}^3 \exp\left(-{\frac{\norm{\boldsymbol x- \boldsymbol x_i}}{\sigma^2}}\right),
\end{equation}
where $\sigma = 1/3$ and source points $\boldsymbol x_1 = (0, 0, 0)$, $\boldsymbol x_2 = (0.25, 0.85, 0.85)$ and $\boldsymbol x_3 = (0.6, 0.4, 0.4)$. In two dimensions, the source points are projected onto the $xy$-plane at $z=0$.

Our results are summarized in \Cref{tab:cubecp,tab:cubemvp}, where for each polynomial degree convergence steps for discretizations on mesh level $L$ are shown with $10^6$ to $10^9$ degrees of freedom in two dimensions, with $10^5$ to $10^9$ in three dimensions. First, we observe that all step counts are independent of the mesh level. Thus, we confirm that we have uniform convergence with respect to mesh size.
The additive cell based smoother (ACS) requires a relaxation parameter $\omega = 0.7$. \Cref{tab:cubecp} shows a slight growth of the number of iteration steps with polynomial degree. It takes about twice as many steps as the multiplicative version (MCS) with $\omega = 1$. Given that MCS with red-black coloring in~\eqref{eq:coloredschwarz} needs two applications of the operator $A_\ell$ in each step, the two cell based smoothers compare at similar levels.

\begin{table}[tp]
\caption{Fractional iterations $\fracsteps$ for cell based smoothers. Multigrid preconditioner with additive smoother (ACS) for CG solver and with multiplicative smoother (MCS) for GMRES solver with relative accuracy of $10^{-8}$. Entries ``---'' not computed.}
\label{tab:cubecp}
\centering
\begin{tabular}{c ccc}
\toprule 
Level $L$ & \multicolumn{3}{c}{Iterations $\fracsteps$ (ACS)} \\
\cmidrule(lr){1-1} \cmidrule(lr){2-4}
2D & $k=3$ & $k=7$ & $k=15$ \\
\midrule 
\myrowcolor
6 & ---  & ---  & 25.4 \\
7 & ---  & 18.7 & 25.4  \\
\myrowcolor
8 & 14.5 & 18.7 & 25.4  \\
9 & 14.5 & 18.7  & 25.4 \\
\myrowcolor
10 & 14.4 & 18.7  & ---  \\
11 & 14.4 & ---  & ---  \\
\midrule 
3D & $k=3$ & $k=7$ & $k=15$\\
\midrule 
\myrowcolor
1 & ---  & ---  & 28.5  \\
2 & ---  & 21.9 & 29.5  \\
\myrowcolor
3 & 17.1 & 21.9 & 29.4  \\
4 & 17.2 & 22.3 & 29.5 \\
\myrowcolor
5 & 17.2 & 22.3  & ---  \\
6 & 17.1 & ---   & ---  \\
\bottomrule
\end{tabular}
\label{tab:cubeacp}
\hfill
\begin{tabular}{c ccc}
\toprule 
Level $L$ & \multicolumn{3}{c}{Iterations $\fracsteps$ (MCS)} \\
\cmidrule(lr){1-1} \cmidrule(lr){2-4}
2D & $k=3$ & $k=7$ & $k=15$ \\
\midrule 
\myrowcolor
6 &  --- &  --- & 12.6 \\
7 &  --- &  9.7 & 12.5 \\
\myrowcolor
8 &  7.3 &  9.6 & 12.5  \\
9 &  7.3 & 9.6  & 12.4  \\
\myrowcolor
10 &  7.2 & 9.5  & ---  \\
11 &  6.9 & ---  & ---  \\
\midrule 
3D & $k=3$ & $k=7$ & $k=15$ \\
\midrule 
\myrowcolor
1 & ---  & ---  & 15.7 \\
2 & ---  & 11.8 & 15.9 \\
\myrowcolor
3 &  8.6 & 11.8 & 15.8 \\
4 &  8.6 & 11.7 & 15.7 \\
\myrowcolor
5 &  8.6 & 11.6  & ---  \\
6 &  8.5 & ---  & ---  \\
\bottomrule
\end{tabular}
\end{table}

Next, we consider the vertex patch smoother. It is well known that the relaxation parameter for additive methods with overlap has to be chosen smaller than $2^{-d}$, which slows down convergence considerably. Therefore, we only consider the multiplicative version (MVS) here. We use the regular coloring in \Cref{fig:multi-red-black} to process as many patches in parallel as possible. We also compare to the graph-based coloring from~\cite{turcksin16}.
As pointed out in~\Cref{sec:schwarz}, the reordering of local solvers in the multiplicative method, in~\Cref{tab:cubemvp} attributed to different coloring algorithms, affects the smoothing. Both coloring schemes yield iteration counts close to two, with a slight advantage for the algorithm with less colors. This is almost a direct solver. Thus, we conclude that in particular the multiplicative vertex patch smoother (MVS) is a mathematically very well suited algorithm if we manage to implement it efficiently. However, the computational effort of one smoothing step is quite high compared to other smoothers, such that we must compare the total effort to decide on the optimal version. 

\begin{table}[tp]
\caption{Fractional iterations $\fracsteps$ for multiplicative vertex patch smoother (MVS). GMRES solver with relative accuracy $10^{-8}$ preconditioned by multigrid with MVS. Entries ``---'' not computed and ``Colors'' refers to the mesh on level $L$.}
\label{tab:cubemvp}
\centering
\begin{subtable}[t]{0.49\textwidth}
\begin{adjustbox}{max width=\textwidth}
\begin{tabular}{c cccc}
\toprule 
Level $L$ & \multicolumn{3}{c}{Convergence steps $\fracsteps$} & Colors \\
\cmidrule(lr){1-1} \cmidrule(lr){2-4} \cmidrule(lr){5-5}
2D & $k=3$ & $k=7$ & $k=15$ \\
\midrule 
\myrowcolor
6 & --- & --- & 1.7 & 8\\
7 & --- & 2.1 & 1.7 & 8\\
\myrowcolor
8 & 2.5 & 2.1 & 1.7  & 8\\
9 & 2.5 & 2.1  & 1.7  & 8\\
\myrowcolor
10 & 2.4 & 2.0  & ---  & 8\\
11 & 2.4 & ---  & ---  & 8\\
\midrule 
3D & $k=3$ & $k=7$ & $k=15$ \\
\midrule 
\myrowcolor
1 & ---  & ---  & 1.5 & 15 \\
2 & ---  & 2.0 & 1.7 & 16 \\
\myrowcolor
3 & 2.4 & 2.1 & 1.7 & 16 \\
4 & 2.4 & 2.1 & 1.7 & 16 \\
\myrowcolor
5 & 2.4 & 2.1 & --- & 16 \\
6 & 2.4 & --- & --- & 16 \\
\bottomrule
\end{tabular}
\end{adjustbox}
\caption{ MVS based on a minimal coloring.}
\label{tab:cubemvpminimal}
\end{subtable}
\begin{subtable}[t]{0.49\textwidth}
\begin{adjustbox}{max width=\textwidth}
\begin{tabular}{c cccc}
\toprule 
Level $L$ & \multicolumn{3}{c}{Convergence steps $\fracsteps$} & Colors \\
\cmidrule(lr){1-1} \cmidrule(lr){2-4} \cmidrule(lr){5-5}
2D & $k=3$ & $k=7$ & $k=15$ \\
\midrule 
\myrowcolor
\myrowcolor
6 & --- & --- & 2.4 & 17\\
7 & --- & 2.6 & 2.4 & 17\\
\myrowcolor
8 & 2.9 & 2.6 & 2.4  & 17\\
9 & 2.8 & 2.5  & 2.3  & 17 \\
\myrowcolor
10 & 2.8 & 2.5  & ---  & 17 \\
11 & 2.9 & ---  & ---  & 17 \\
\midrule 
3D & $k=3$ & $k=7$ & $k=15$ \\
\midrule 
\myrowcolor
1 & ---  & ---  & 1.7 & 19 \\
2 & ---  & 2.3 & 2.1 & 36 \\
\myrowcolor
3 & 2.7 & 2.4 & 2.2 & 50 \\
4 & 2.7 & 2.5 & 2.2 & 52 \\
\myrowcolor
5 & 2.7 & 2.5  & --- & 53 \\
6 & 2.7 & ---  & --- & 53 \\
\bottomrule
\end{tabular}
\end{adjustbox}
\caption{MVS based on coloring by DSATUR.}
\label{tab:cubemvpgraph}
\end{subtable}
\end{table}

\section{Tensor Product Elements}
\label{sec:tpelements}
In recent years, the structure of tensor product polynomials on quadrilateral and hexahedral cells and their evaluation and integration by sum factorization has been exploited in the development of highly efficient codes for modern hardware, see for instance~\cite{KronbichlerKormann12,MuethingPiatkowskiBastian17,vos10}.

Let $\poly k$ be the space of polynomials in one variable of degree up to $k$ with basis $\{ \hat\phi_i \}_{i=0,\ldots,k}$. We use the Lagrange basis in Gauss-Lobatto support points for stability.
Then, we define the tensor product polynomial space
\begin{equation}\label{eq:tppoly}
    \tppoly{k} = \bigotimes_{\tau=1}^d \poly{k},
\end{equation}
with its basis
\begin{equation}\label{eq:basisref}
  \hat\varphi_{\mulindex{i}}(\bfx)
  := \hat\phi_{i_1}\otimes\dots\otimes\hat\phi_{i_d}(\bfx)
  = \prod_{\tau=1}^d \hat\phi_{i_\tau}(x_\tau), \quad i=0,\ldots,N-1.
\end{equation}
where $N-1=(k+1)^d$. Here, we have adopted a multi-index notation, such that each basis function is characterized by $d$ indices. These can be unrolled into a linear index $i \in \{0,\dots,N-1\}$, for instance by the lexicographic mapping 
\begin{equation}\label{eq:multiindex}
  i = \sum_{\tau=1}^d i_\tau (k+1)^{\tau-1}.
\end{equation}
This defines a polynomial shape function space $V(\refcell)$ on the reference cell $\refcell$. The polynomial shape function spaces $V(\cell)$ and its basis $\{\varphi_{\cell,i}\}$ on the mesh cell $\cell$ are obtained by composition with the cell mapping $\mapping\cell$, that is, $\varphi_{\cell,i} (\bfx) := \hat\varphi_i \circ \mapping{\cell}^{-1} (\bfx)$.
Similarly, we define a $d$-dimensional quadrature rule on $\refcell$ as the $d$-fold tensor product of a one-dimensional on the unit interval. If the one-dimensional quadrature formula has abscissas and weights $\bigl\{(\xref_{q_\tau},\qweight_{q_\tau})\bigr\}$ on the interval $[0,1]$, we let
\begin{gather}
  \bfxref_{\multiindex{q}} = (\hat{x}_{q_1},\dots,\hat{x}_{q_d})^T,
  \qquad
  \qweight_{\multiindex{q}} = \prod_{\tau=1}^d \qweight_{q_\tau}.
\end{gather}

\subsection{Operator application by sum factorization}
\label{sec:matrixfree}

A matrix-free finite element implementation with sum factorization is easiest discussed using the mass matrix as example. Instead of assembling the mass matrix $M_l$ on the mesh on level $\ell$, integration and application to a vector in coefficient space are folded into one operation. Let $\doftransfer$ denote the transfer from global degrees of freedom on $\tria \ell$ to local degrees of freedom on the cell $\cell$. Then,
\begin{equation}
  M_\ell u_\ell = \sum_{\cell\in\tria\ell} \doftransfer^T M_{\ell,\cell} \doftransfer u_\ell = \sum_{\cell\in\tria\ell} \doftransfer^T M_{\ell,\cell} u_{\cell}.
\end{equation}
Here, $u_\ell$ denotes the coefficient vector of the finite element function $u_\ell$ with respect to the chosen basis.
The restriction $u_{\cell} \in V(K)$ of the finite element function $u_\ell$ to the cell $\cell$ is determined by real-valued coefficients $u_\cell$, reshaped as order-$d$ tensor $U_{\mulindex i}$, where here and below we suppress the cell index $\cell$.
\begin{equation}\label{eq:dofexpansion}
  u_{\cell}(\bfx)
  = \sum_{\mulindex i=0}^{k} U_{\mulindex i} \varphi_{\cell,\mulindex i} (\bfx)
  = \sum_{\mulindex i=0}^{k} U_{\mulindex i} \hat{\varphi}_{\mulindex i} \circ \mapping{K}^{-1} (\bfx).
\end{equation}
We assume the number of univariate quadrature points to be almost identical with the polynomial degree $k$. Understanding the local finite element interpolation $u_\cell$ as pullback by the mapping $\mapping{\cell}$, the evaluation of $u_{\cell}$ in the standard form~\eqref{eq:dofexpansion} in all quadrature points requires $\bigo (k^{2d})$ arithmetic operations.

Exploiting the tensor product form of basis functions and quadrature formula, we can reduce the complexity to $\bigo (d k^{d+1})$ by means of \emph{sum factorization}. This technique was first introduced in the spectral element community in~\cite{orszag80} and later extended to DG methods, see for instance \cite{vos10}. By factorizing common indices in~\eqref{eq:dofexpansion} along each dimension we obtain $d$ one-dimensional interpolations
\begin{align}\label{eq:dofexpansionsf}
  \begin{aligned}
    \hat{U}_{q_1,\ldots,q_d}
    = u_{\cell} (\bfx_{\mulindex q}) 
    = \sum_{i_d=0}^{k} \hat\phi_{i_d} (\xref_{q_d}) \cdots
      \sum_{i_2=0}^{k} \hat\phi_{i_2} (\xref_{q_2})
      \sum_{i_1=0}^{k} U_{i_1, \ldots, i_d} \hat\phi_{i_1} (\xref_{q_1})
  \end{aligned}
\end{align}
where implicitly $\bfxref_q = \mapping\cell^{-1} (\bfx_q)$ is used. The order-$d$ tensor $\hat{U}$ is successively obtained by computing the $d$ sum-factors from right to left. In other words, each sum-factor results in an intermediate tensor with degree of freedom and quadrature indices mixed obtained by the contraction of the previous order-$d$ tensor and the matrix composed of $\hat\phi_{i_\tau} (\xref_{q_\tau})$ for all $i_\tau,q_\tau$. A change of variables with the cell mapping $\mapping\cell$ and application of the quadrature formula yields
\begin{align}
\label{eq:massmatvec}
\begin{aligned}
  \left( M_\cell u_\cell \right)_{\mulindex i}
  &= \int_{\cell} \varphi_{K,\mulindex i} u_{\cell} \dbfx\\
  &= \sum_q \hat\varphi_{\mulindex i}(\bfxref_q) u_{\cell}(\bfx_q) \determinant\jacobianr_\cell(\bfxref_q) \qweight_q,
\end{aligned}
\end{align}
where we use that $\nabla F$ is chosen with positive determinant.
Switching to the multi-indices of $q$, we introduce the order-$d$ mapping tensor $\tensorJxW$ with components
\begin{gather}
    \tensorJxW_{\mulindex q} = \det\jacobianr_\cell(\xref_{q_1},\ldots,\xref_{q_d}) \qweight_{q_1,\ldots,q_d}.
\end{gather}
This time, factorizing the sum over quadrature points along each dimension and using~\eqref{eq:dofexpansionsf} we transform the integration against all test functions $\varphi_{\cell,i}$ in~\eqref{eq:massmatvec}
\begin{equation}\label{eq:massmatvecsf}
  \left( M_\cell u_\cell \right)_{\mulindex i}
  =
  \sum_{q_d} \hat\phi_{i_d} (\xref_{q_d})
  \cdots
  \sum_{q_2} \hat\phi_{i_2} (\xref_{q_2})
  \sum_{q_1} \bigl(\tensorJxW \hadap \hat{U}\bigr)_{q_1,\ldots,q_d} \hat\phi_{i_1}(\xref_{q_1})
\end{equation}
with $\hadap$ being the entrywise product (also known as the Hadamard product). In total, the local matrix-vector multiplication with $M_\cell$ is performed at the cost of $\bigo (dk^{d+1})$ arithmetic operations. In this context, we refer to $M_l$ as a matrix-free finite element operator.

A similar expression can be derived for the Laplacian or a general second order elliptic operator on arbitrary quadrilaterals and hexahedra. The formula becomes more complicated then since it involves a matrix-valued mapping tensor $T$ and a vector-valued interpolation tensor $\hat{U}$ due to the gradients of ansatz and test functions.

\subsection{Fast diagonalization}
\label{sec:fastdiagsmoother}

In \Cref{sec:robustcart}, we have presented robust Schwarz smoothers with low iteration counts. However, the naive computation of local inverses $A_{j}^{-1}$ requires $\bigo (k^{3d})$ arithmetic operations, while matrix-free operator application $A_{j} v_{j}$ costs only $\bigo (dk^{d+1})$. Therefore, explicit inversion should be avoided as well as multiplication with an inverse with $\bigo (k^{2d})$ operations.

The \emph{fast diagonalization method} introduced in \cite{lynch64} is an efficient inversion algorithm for matrices $A$ with a rank-d Kronecker decomposition of the form
\begin{equation}
    \label{eq:KDseparable}
    A =
    M^{(d)} \otimes \cdots \otimes M^{(2)} \otimes A^{(1)}
    + \ldots 
    + A^{(d)} \otimes M^{(d-1)} \otimes \cdots \otimes M^{(1)}.
\end{equation}
Assuming the matrices $M^{(\tau)}$ are symmetric, positive definite and the matrices $A^{(\tau)}$ are symmetric, the generalized eigenvalue problems
\begin{equation}
    \label{eq:generalizedEVP}
    (Z^{(\tau)})^T A^{(\tau)} Z^{(\tau)} = \Lambda^{(\tau)},
    \quad
    (Z^{(\tau)})^T M^{(\tau)} Z^{(\tau)} = I^{(\tau)},
    \quad
    \tau = 1,\ldots,d,
\end{equation}
 are well-defined. Here $Z^{(\tau)}$ is the orthogonal matrix of generalized eigenvectors and $I^{(\tau)}$ is the identity matrix of appropriate size.
Using the mixed-product property of the Kronecker product we see that the Kronecker product of the generalized eigenvectors, namely $Z := Z^{(d)} \otimes\dots\otimes Z^{(1)}$, are the eigenvectors of $A$
\begin{equation}\label{eq:localsolverEVP}
    Z^T A Z
    = I^{(d)} \otimes\dots\otimes I^{(2)} \otimes \Lambda^{(1)}
    + \ldots
    + \Lambda^{(d)} \otimes I^{(d-1)} \otimes\dots\otimes I^{(1)}
    =: \Lambda,
\end{equation}
such that the inverse of $A$ is
\begin{equation}
    \label{eq:fastinverse}
     A^{-1} = Z \Lambda^{-1} Z^T.
\end{equation}
The rank-$1$ Kronecker decomposition $Z^{(d)} \otimes \cdots \otimes Z^{(1)}$ plays a key role in obtaining a fast inversion algorithm. First, the assembly and inversion of $A$ boils down to the assembly and subsequent computation of the generalized eigendecomposition of $d$ one-dimensional problems, respectively. Second, we only store $d$ one-dimensional eigenvector matrices $Z^{(\tau)}$ and eigenvalues $\Lambda^{(\tau)}$. Then, the matrix-vector multiplication $Zu$ profits from its Kronecker decomposition in terms of sum factorization
\begin{equation}\label{eq:eigenvectormatvec}
  (Z^{(d)} \otimes \cdots \otimes Z^{(1)}u)_m =
  \sum_{n_d} Z^{(d)}_{m_d,n_d}
  \cdots
  \sum_{n_2} Z^{(2)}_{m_2,n_2}
  \sum_{n_1} Z^{(1)}_{m_1,n_1} \tensoru_{n_1,\ldots,n_d} ,
\end{equation}
where the order-$d$ tensor $\tensoru$ is the multi-index reshaping of the vector $u$.

The Laplacian is a separable differential operator on rectangles and bricks, but it remains to argue that the discontinuous Galerkin formulation on a cell or a vertex patch is as well. Indeed, this is true on Cartesian meshes only. Let the cell $\cell$ be of the form $\cell = \interval_1 \times \cdots \times \interval_d$ with intervals $\interval_\tau = [x_{\tau,0},x_{\tau,1}]$ of length $h_\tau$.
\begin{equation}\label{eq:bulkcart}
  \int_\cell \grad \varphi_{K,i} \cdot \grad \varphi_{K,\iota} \dbfx
  = \sum_{t=1}^d \left( \int_{\interval_t} \phi_{K,i_t}^\prime \phi_{K,\iota_t}^\prime \dx_t
  \prod_{\tau=1, \tau \neq t}^d \int_{\interval_\tau} \phi_{K,i_\tau} \phi_{K,\iota_\tau} \dx_\tau \right).
\end{equation}
Then, the bulk integral in~\eqref{eq:bfsipg} on cell $\cell$ is the sum of products~\eqref{eq:bulkcart} alternating with the dimension. Each product is factorized as a one-dimensional bulk integral and $d-1$ remaining $L^2$-inner products of one-dimensional shape functions.

In general the shape function gradients $\grad \varphi_{\cell,i} \circ \mapping\cell$ are determined by means of the chain rule $\jacobianr_\cell^{-T} \gradref \hat\varphi_{\cell,i}$. The Jacobian of the Cartesian mapping is the constant, diagonal matrix $\diag \left( h_1,\ldots,h_d \right)$ such that the univariate mass and interior stiffness matrices $M^{(\tau)}$ and $L^{(\tau)}$ are
\begin{equation}\label{eq:1Dstiffnessmass}
  \left( M^{(\tau)} \right)_{\iota,i}
  = \sum_{q_\tau} \hat{\phi}_i (\xref_{q_\tau}) \hat{\phi}_\iota (\xref_{q_\tau}) h_\tau \qweight_{q_\tau},
  \quad
  \left( L^{(\tau)} \right)_{\iota,i}
  = \sum_{q_\tau} \frac{1}{h_\tau} \hat{\phi}_i^\prime (\xref_{q_\tau}) \hat{\phi}_\iota^\prime  (\xref_{q_\tau}) \qweight_{q_\tau},
\end{equation}
respectively, where the quadrature rule is defined on the reference interval $[0,1]$. 
The two faces $\face_p$ of $\cell$ associated to dimension $\tau$ are as well a Cartesian product
\begin{equation}
  \face_p :=
  \interval_1 \times \cdots \times \interval_{\tau-1} \times \{ x_{\tau,p}\} \times \interval_{\tau+1} \times \cdots \times \interval_d ,
  \quad p = 0,1 ,
\end{equation}
such that the face normals are constant and aligned with coordinate direction. We omit the subscript $p$ of the face $\face$ when it is clear from the context. We obtain a similar splitting for the consistency, adjoint consistency and penalty integrals in~\eqref{eq:bfsipg}. The univariate consistency and point mass matrices $G_{\face}^{(\tau)}$ and $M^{(\tau)}_{\face}$ for the faces orthogonal to coordinate direction $\tau$ are obtained by
\begin{equation}\label{eq:1Dconsistencypenalty}
  \left( G_{\face,p}^{(\tau)} \right)_{\iota,i}
  = (-1)^{p+1} \frac{\eta_\face}{h_\tau} \hat{\phi}_i^\prime (p) \hat{\phi}_\iota (p) ,
  \quad
  \left( M_{p}^{(\tau)} \right)_{\iota,i}
  = \hat{\phi}_i (p) \hat{\phi}_\iota (p) ,
  \quad
  p = 0,1 ,
\end{equation}
where $\eta_\face=1$ on any face $\face$ at the physical boundary and $\eta_\face=\nicefrac{1}{2}$ otherwise. The Nitsche contributions are summed up to obtain
\begin{equation}\label{eq:1Dnitsche}
    N_{\face,p}^{(\tau)} = \ip_\face M^{(\tau)}_{p} - G_{\face,p}^{(\tau)} - \bigr(G_{\face,p}^{(\tau)}\bigl)^T, \quad p = 0,1.
\end{equation}
Hence, the local solvers on cells admit a Kronecker decomposition of the form~\eqref{eq:KDseparable} with
\begin{equation}
    \label{eq:1Dsipg}
    A^{(\tau)}
    = L^{(\tau)}
    + N_{\face,0}^{(\tau)} + N_{\face,1}^{(\tau)}.
\end{equation}

Cartesian vertex patches are determined by the Cartesian product of intervals $\interval_\tau$. Each interval is defined by the disjoint union $\interval_{\tau,+} \dot\cup \interval_{\tau,-}$ with subintervals $\interval_{\tau,\pm} := [a_{\tau,\pm},b_{\tau,\pm}]$ of length $h_{\tau,\pm}$ and $b_{\tau,+} = a_{\tau,-}$. Besides the interior contributions~\eqref{eq:1Dsipg} on both subintervals, denoted as $A_+^{(\tau)}$ and $A_-^{(\tau)}$, the contributions from the interface between $\interval_{\tau,+}$ and $\interval_{\tau,-}$ have to be considered. The univariate consistency and point mass matrices $N_{\face,\pm\mp}^{(\tau)}$ and $M_{\pm\mp}^{(\tau)}$, respectively, are given by
\begin{align}
  \label{eq:1Dconsistencyinterface}
  \left( G_{\face,+-}^{(\tau)} \right)_{\iota,i}
  &= \frac{1}{2h_{\tau,-}} \bigl( \hat{\phi}_i \bigr)^\prime (0) \hat{\phi}_\iota (1) ,
  &
    \left( G_{\face,-+}^{(\tau)} \right)_{\iota,i}
  &= - \frac{1}{2h_{\tau,+}} \bigl( \hat{\phi}_i \bigr)^\prime (1) \hat{\phi}_\iota (0) ,
  \\ \label{eq:1Dpenaltyinterface}
  \left( M_{+-}^{(\tau)} \right)_{\iota,i}
  &= \hat{\phi}_i (0) \hat{\phi}_\iota (1) ,
  &
    \left( M_{e,-+}^{(\tau)} \right)_{\iota,i}
  &= \hat{\phi}_i (1) \hat{\phi}_\iota (0) ,
\end{align}
and the Nitsche terms at the interface are summed up as
\begin{equation}
    N_{\face,\ast}^{(\tau)}
    = \ip_\face M^{(\tau)}_{e,\ast} - G_{e,\ast}^{(\tau)} - \bigl(G_{e,\ast}^{(\tau)}\bigr)^T,
\end{equation}
replacing $\ast$ by $+-$ or $-+$, respectively. Therefore, the local solvers on vertex patches admit a Kronecker decomposition of the form~\eqref{eq:KDseparable} with
\begin{equation}
    \label{eq:1DsipgVP}
    A^{(\tau)} = 
    \begin{bmatrix}
    A_{+}^{(\tau)} &  A_{+-}^{(\tau)} \\
    A_{-+}^{(\tau)} & A_{-}^{(\tau)} \\
    \end{bmatrix},
    \qquad
    M^{(\tau)} = 
    \begin{bmatrix}
    M_{+}^{(\tau)} & 0 \\ 
    0 & M_{-}^{(\tau)} \\
    \end{bmatrix},
\end{equation}
where $M_{\pm}^{(\tau)}$ are the mass matrices \eqref{eq:1Dstiffnessmass} and $A_{\pm}^{(\tau)}$ the stiffness matrices \eqref{eq:1Dsipg} on $\interval_{\tau,\pm}$, respectively. The interior penalty interface matrices $A_{\pm\mp}^{(\tau)}$ are defined by
\begin{equation}\label{eq:1Dsipginterface}
  A_{\ast}^{(\tau)}
  =  N_{e,\ast}^{(\tau)} +  N_{e,\ast}^{(\tau)} ,
\end{equation}
replacing $\ast$ by $+-$ or $-+$, respectively.

\subsection{Computational effort}
\label{sec:compeffort}
We compare the computational effort of the fast tensor product smoothers with the other components of the multigrid scheme. The experimental setup is the same as in \Cref{sec:robustcart}. Consequently, the operation counts presented here are consistent with the iteration counts there. We compute on a three-dimensional mesh $\tria 3$ with $8^4=4096$ cells obtained from a single vertex patch (coarse grid) by three consecutive global refinements.
We use implementations based on \texttt{deal.II} \cite{dealII91} and in particular its \texttt{MatrixFree} framework.

\begin{table}[tp]
\centering
\caption{Asymptotic work load per cell (ACS and residual) or vertex patch (AVS) of additive smoothers in three dimensions. Leading order of the setup is independent of the dimension, consisting of $d$ one-dimensional eigenvalue solvers.}
\label{tab:cmplxconstant}
\begin{tabular}{l rrrrrrr l}
\toprule 
Method & \multicolumn{7}{c}{Factor $\ccmplx$} & Order \\
\cmidrule(lr){1-1} \cmidrule(lr){2-8} \cmidrule(lr){9-9}
Degree $k$: & \multicolumn{1}{c}{7} & \multicolumn{1}{c}{11} & \multicolumn{1}{c}{15} & \multicolumn{1}{c}{19} & \multicolumn{1}{c}{23} & \multicolumn{1}{c}{27} & \multicolumn{1}{c}{31} \\
\midrule 
$A_\ell u_\ell$ & 32 & 25 & 21 & 19 & 18 & 18 & 17 & $\times k^{d+1}$ \\
\myrowcolor
ACS: $S_{\ell}(u_{\ell},b_{\ell})$ & 45 & 38 & 34 & 32 & 31 & 30 & 29 & $\times k^{d+1}$ \\
\myrowcolor
\hphantom{ACS: }--- local solvers & 12 & 12 & 12 & 12 & 12 & 12 & 12 & $\times k^{d+1}$ \\
\myrowcolor
\hphantom{ACS: }setup of $S_\ell$ & 96 & 77 & 67 & 63 & 57 & 54 & 52 & $\times k^3$ \\
AVS: $S_{\ell}(u_{\ell},b_{\ell})$ & 236 & 227 & 222 & 219 & 218 & 217 & 216 & $\times k^{d+1}$ \\
\hphantom{AVS: }--- local solvers & 195 & 195 & 195 & 195 & 195 & 195 & 195 & $\times k^{d+1}$ \\
\hphantom{AVS: }setup of $S_\ell$ & 506 & 419 & 369 & 393 & 365 & 354 & 341 & $\times k^3$ \\
\bottomrule
\end{tabular}
\end{table}

First, we confirm the asymptotic complexity of the fast tensor product smoothers with respect to polynomial degree.
The number of floating point operations $n_{\textit{FLOP}}$ is determined by means of the performance monitoring tool \texttt{likwid-perfctr} \cite{treibig10}.
In \Cref{tab:cmplxconstant}, we report the factors
\begin{equation}
    \label{eq:cmplxconstant}
    \ccmplx = \frac{n_{\textit{FLOP}}}{\normalizationcmplx \times k^{\ordercmplx}},
\end{equation}
which are obtained from normalizing FLOP counts by the number $\normalizationcmplx$ of subdomains (cells or vertex patches) and the expected complexity $k^{\ordercmplx}$ in the polynomial degree.

In \Cref{tab:cmplxconstant}, we confirm the asymptotic behavior (last column) of the additive Schwarz smoothers in \Cref{alg:addsmoother} over a wide range of polynomial degrees. The arithmetic effort for the smoother consists of two parts, namely the one-time setup cost and the smoothing operation $S_\ell(u_\ell,b_\ell) $  in each step. 
The setup consists mainly of $d$ one-dimensional eigenvalue problems including integration and solving (\texttt{LAPACK} routine DSYGV) with order $k^3$ operations. The integration cost decreases for higher order polynomials since the contribution of lower order face terms becomes less prominent.
The effort of one additive smoothing step $S_{\ell}(u_{\ell},b_{\ell})$ is determined by the cost of applying all local solvers and updating the residual, that is a single operator application $A_\ell u_\ell$.
The local solvers scale strictly with $k^{d+1}$, while the normalized numerical effort for residuals decreases with increasing polynomial degree as the face integrals of the DG discretization lose weight. Thus, compared to a operator application, the additive smoothing step becomes cheaper with increasing degree.
Considering now the vertex patch, the degrees of freedom in each dimension are doubled. Therefore, applying the local solvers of AVS is computationally $16$ times more expensive than in ACS. Both the matrix-free operator application and the smoothing step crucially benefit from the tensor structure. 

\begin{table}[tp]
\centering
\caption{Arithmetic operations, additive smoothers $S_{\ell;\textit{ad}}$ vs. operator application $A_\ell$ in MFLOPS, $\ell=3$. sACS is the standard smoother without exploiting tensor structure.}
\label{tab:flop-additive}
\begin{adjustbox}{max width=\textwidth}
\begin{tabular}{l rr rrr rrr}
\toprule 
Method & \multicolumn{2}{c}{sACS} & \multicolumn{3}{c}{ACS} & \multicolumn{3}{c}{AVS} \\
\cmidrule(lr){1-1} \cmidrule(lr){2-3} \cmidrule(lr){4-6} \cmidrule(lr){7-9}
Degree $k$: & \multicolumn{1}{c}{3} & \multicolumn{1}{c}{7} & \multicolumn{1}{c}{3} & \multicolumn{1}{c}{7} & \multicolumn{1}{c}{15} & \multicolumn{1}{c}{3} & \multicolumn{1}{c}{7} & \multicolumn{1}{c}{15} \\
\midrule 
$A_\ell u_\ell$ & 59 & 545 & 59 & 545 & 5,819 & 59 & 545 & 5,819 \\
\myrowcolor
$S_{\ell;\textit{ad}}(u_{\ell},b_{\ell})$ & 129 & 5,474 & 74 & 763 & 9,176 & 229 & 3,235 & 48,629 \\
--- local solvers & 69 & 4,861 & 13 & 205 & 3,256 & 168 & 2,677 & 42,709 \\
\myrowcolor
--- residual & 60 & 551 & 60 & 551 & 5,869 & 60 & 551 & 5,869 \\
setup of $S_{\ell;\textit{ad}}$ & 176,591 & 52,073,595 & 37 & 206 & 1,143 & 163 & 866 & 5,013 \\
\bottomrule
\end{tabular}
\end{adjustbox}
\end{table}

In \Cref{tab:flop-additive}, we compare the number of arithmetic operations for constituents of the additive smoothing and, in particular, we compare to the non-tensorized cell based smoother (sACS). Reading columns two to five, the benefits of the fast tensor product smoothers are exposed. The number of operations to setup the smoother sACS are 3000 times higher than a single operator application for tricubic shape functions, even 95000 times higher for $\tppoly{7}$. Clearly, non-tensorized Schwarz smoothers are infeasible since they obliterate the advantages of matrix-free methods. As columns four and five show, the setup cost of the fast tensor product cell based smoother (ACS) is already less than a single matrix application and therefore almost negligible.
We see that one smoothing step needs about $\nicefrac54$ of the number of operations of a matrix application for tricubic shape functions, about $\nicefrac85$ for $\tppoly{15}$, both well bounded below 2.
The cost for the smoother in each step is dominated by the computation of the residual. 
When we compare to the additive vertex patch, we realize that applying the local solvers costs about 13 times as much. As already mentioned a single vertex patch solver is $2^{d+1}$ times more expensive due to the doubled size of the patch in each direction. Hence, we would expect to see a factor 16 if the number of cells and patches would be the same, but $\tria 3$ has 3375 vertex patches and 4096 cells and $16 \times \nicefrac{3375}{4096} \approx 13$.

\begin{table}[tp]
\centering
\caption{Arithmetic operations, multiplicative smoothers $S_{\ell;\textit{mu}}$ vs. operator application $A_\ell$ in MFLOPS, $\ell=3$.}
\label{tab:flop-multiplicative}
\begin{adjustbox}{max width=\textwidth}
\begin{tabular}{lrrrrrr}
\toprule 
Method & \multicolumn{3}{c}{MCS} & \multicolumn{3}{c}{MVS} \\
\cmidrule(lr){1-1} \cmidrule(lr){2-4} \cmidrule(lr){5-7}
Degree $k$: & \multicolumn{1}{c}{3} & \multicolumn{1}{c}{7} & \multicolumn{1}{c}{15}
& \multicolumn{1}{c}{3} & \multicolumn{1}{c}{7} & \multicolumn{1}{c}{15} \\
\midrule 
$A_\ell u_\ell$ & 59 & 545 & 5819 & 59 & 545 & 5819 \\
\myrowcolor
$S_{\ell;\textit{mu}}(u_{\ell},b_{\ell})$ & 135 & 1,319 & 15,077 & 1,142 & 11,633 & 138,576 \\
--- local solvers & 13 & 209 & 3,287 & 176 & 2,802 & 44,612 \\
\myrowcolor
--- residual & 120 & 1,103 & 11,739 & 965 & 8,825 & 93,913 \\
setup of $S_{\ell;\textit{mu}}$ & 37 & 206 & 1143 & 170 & 903 & 5226 \\
\bottomrule
\end{tabular}
\end{adjustbox}
\end{table}

In \Cref{tab:flop-multiplicative}, we show the same data for multiplicative smoothers. First, we notice that while the cost for a single operator application remains the same, residuals are computed once for each color inside the smoother. Hence, the computation of residuals costs two times more for the cell based smoother (MCS), 16 times more for the vertex patch smoother (MVS). Comparing to the additive counterparts, the number of local solvers is the same, thus, their cost is similar. Further reducing the cost of the multiplicative methods by avoiding the intermediate residual computations is possible, but requires a major change of implementation. 

The significant difference in the computational effort needed for a single smoothing step motivates the discussion of the trade-off between the effort of a smoother and its iteration counts. Thus, we end this subsection on the computational effort in FLOPS with \Cref{tab:full-additive,tab:full-multiplicative}, where we compare the effort for a single operator application to the whole multigrid solver including smoothing and grid transfer.
\begin{table}[tp]
\centering
\caption{Arithmetic operations (MFLOPS) for single V-cycle $\mathrm{MG}_{\ell}$ with additive smoothers $S_{\ell;\textit{ad}}$ and pCG solver to relative accuracy $10^{-8}$.}
\label{tab:full-additive}
\begin{adjustbox}{max width=\textwidth}
\begin{tabular}{lrrrrrr}
\toprule 
Method & \multicolumn{3}{c}{ACS} & \multicolumn{3}{c}{AVS} \\
\cmidrule(lr){1-1} \cmidrule(lr){2-4} \cmidrule(lr){5-7}
Degree $k$: & \multicolumn{1}{c}{3} & \multicolumn{1}{c}{7} & \multicolumn{1}{c}{15} & \multicolumn{1}{c}{3} & \multicolumn{1}{c}{7} & \multicolumn{1}{c}{15} \\
\midrule 
$A_\ell u_\ell$ & 59 & 545 & 5819 & 59 & 545 & 5819 \\
\myrowcolor
$S_{\ell;\textit{ad}}(u_{\ell},b_{\ell})$ & 74 & 763 & 9,176 & 229 & 3,235 & 48,629 \\
$\mathrm{MG}_{\ell}(u_{\ell},b_{\ell})$ & 255 & 2,696 & 34,332 & 600 & 8,194 & 122,110 \\
\myrowcolor
solver & 5,645  & 69,707 & 1,141,043 & 19,763 & 294,632 & 4,780,256 \\
\myrowcolor
--- preconditioner & 4,512 & 57,104 & 959,889 & 17,875 & 275,153 & 4,550,847 \\
\bottomrule
\end{tabular}
\end{adjustbox}
\end{table}
\Cref{tab:full-additive} shows that the effort for the whole solver with additive cell based smoother (ACS) ranges from 96 to 200 times the effort for single operator application for polynomial degrees 3 and 15, respectively. In part, this can be explained by the increase in the number of iteration steps from 18 to 30.
\begin{table}[tp]
\centering
\caption{Arithmetic operations (MFLOPS) for single V-cycle $\mathrm{MG}_{\ell}$ with multiplicative smoothers $S_{\ell;\textit{mu}}$ and pGMRES solver to relative accuracy $10^{-8}$.}
\label{tab:full-multiplicative}
\begin{adjustbox}{max width=\textwidth}
\begin{tabular}{lrrrrrr}
\toprule 
Method & \multicolumn{3}{c}{MCS} & \multicolumn{3}{c}{MVS} \\
\cmidrule(lr){1-1} \cmidrule(lr){2-4} \cmidrule(lr){5-7}
Degree $k$: & \multicolumn{1}{c}{3} & \multicolumn{1}{c}{7} & \multicolumn{1}{c}{15} & \multicolumn{1}{c}{3} & \multicolumn{1}{c}{7} & \multicolumn{1}{c}{15} \\
\midrule 
$A_\ell u_\ell$ & 59 & 545 & 5,819 & 59 & 545 & 5,819 \\
\myrowcolor
$S_{\ell;\textit{mu}}(u_{\ell},b_{\ell})$ & 135 & 1,319 & 15,077 & 1,142 & 11,633 & 138,576 \\
$\mathrm{MG}_{\ell}(u_{\ell},b_{\ell})$ & 394 &3,962 & 47,791 & 2,695 & 27,536 & 330,508 \\
\myrowcolor
solver & 4,572 & 58,147 & 883,421 & 11,062 & 112,118 & 1,004,638 \\
\myrowcolor
--- preconditioner & 3,913 & 50,254 & 773,818 & 10,811 & 109,841 & 986,710 \\
\bottomrule
\end{tabular}
\end{adjustbox}
\end{table}
The multiplicative vertex patch (MVS) smoother in \Cref{tab:full-multiplicative} requires between 170 and 200 times the effort of a single operator application for the iterative solution with much less dependence on the polynomial degree.

These numbers suggest that the multigrid preconditioner compares favorably to the unpreconditioned conjugate gradient method if it needs 200 steps. For comparison, between 285 and 1805 steps by an unpreconditioned conjugate gradient were needed for the same discretizations. Furthermore, the discussion in~\cite{KronbichlerLjungkvist19,ArndtFehnKanschatKormannKronbichlerMunchWallWitte19} shows that using an efficient implementation operator application doesn't dominate run time of the conjugate gradient method anymore which is in favor of our method. Most importantly, we obtain a method which is robust with respect to polynomial degree.

\section{Non-Cartesian Meshes}
\label{sec:unstructured}
Exploiting fast diagonalization is intrinsically connected to separability of the differential operator and therefore to a suitable geometry. Thus, it is almost entirely restricted to rectangular meshes. For mesh cells of more general shape, this concept cannot be applied anymore in its original version. We address this issue by replacing the actual grid cell by a rectangular surrogate for smoothing purposes.

From the point of view of Schwarz methods, we replace the local solvers $P_j$ in equation~\eqref{eq:ritzproj} by approximate local solvers $\widetilde{P}_{j;\ell} \colon V_\ell \to V_{j;\ell}$ defined by
\begin{equation}
    \label{eq:schwarzproj}
    \widetilde{a}_{j;\ell} (\widetilde{P}_{j;\ell} u_\ell, v) = a_\ell(u_\ell, v) \quad \forall v \in V_{j;\ell}. 
\end{equation}
Here, $\widetilde{a}_{j;\ell}$ is an approximation to the original form $a_\ell$ locally on the subdomain, which is separable again.
The convergence theory \cite{ToselliWidlund05} for subspace correction methods requires additionally the \emph{local stability} assumption
\begin{equation}
\label{eq:localstab}
a_\ell(u, u) \leq \eta \widetilde{a}_{j;\ell}(u,u),
\qquad u \in V_{j;\ell}, \quad 1 \leq j \leq J.
\end{equation}
for some $\eta>0$, independent of $j$.

The approximation $\widetilde{a}_{j;\ell}$ is obtained by replacing arbitrary mesh cells by Cartesian surrogate cells. Similar ideas were briefly suggested in \cite{couzy95,fischer00} but not implemented. The construction employs the fact that the Laplacian is rotation invariant and works as follows.
\begin{figure}[tbh]
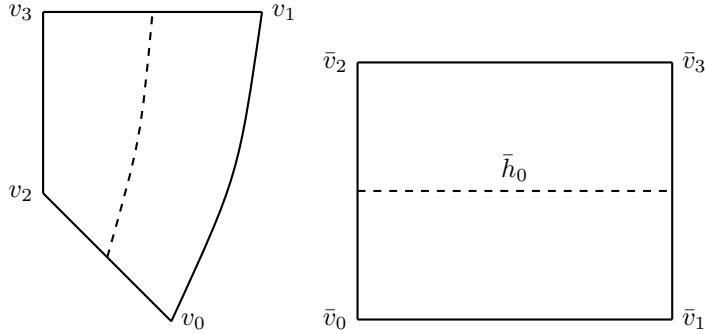

    \centering
    \includegraphics[height=.35\textwidth]{grids/realpatch.tikz}
    \includegraphics[height=.3\textwidth]{grids/surpatch.tikz}
\caption{Illustration of a mesh cell and its surrogate cell.}\label{fig:surrogatecell}
\end{figure}
In \Cref{fig:surrogatecell}, let $\ell_{ij}$ be the arc length of the possibly curved edge connecting the vertices $v_i$ and $v_j$. This arc length can be approximated by a Gauss-Lobatto formula of sufficient order or simply by taking the distance of the two vertices.
Then, define
\begin{gather}
    \label{eq:avglength}
    \bar h_0 = \frac{\ell_{01}+\ell_{23}}{2}
    \quad\text{and}\quad
    \bar h_1 = \frac{\ell_{02}+\ell_{13}}{2}.
\end{gather}
Thus, we have obtained the dimensions of the surrogate rectangle in $x$- and $y$-direction. Its position and orientation are determined by placing $\bar v_0$ and $\bar v_1$ on the $x$-axis and $\bar v_0$ and $\bar v_2$ on the $y$-axis.

\begin{figure}[tbh]
  \centering
  \begin{subfigure}{0.48\textwidth}
  \includegraphics[width=0.48\textwidth]{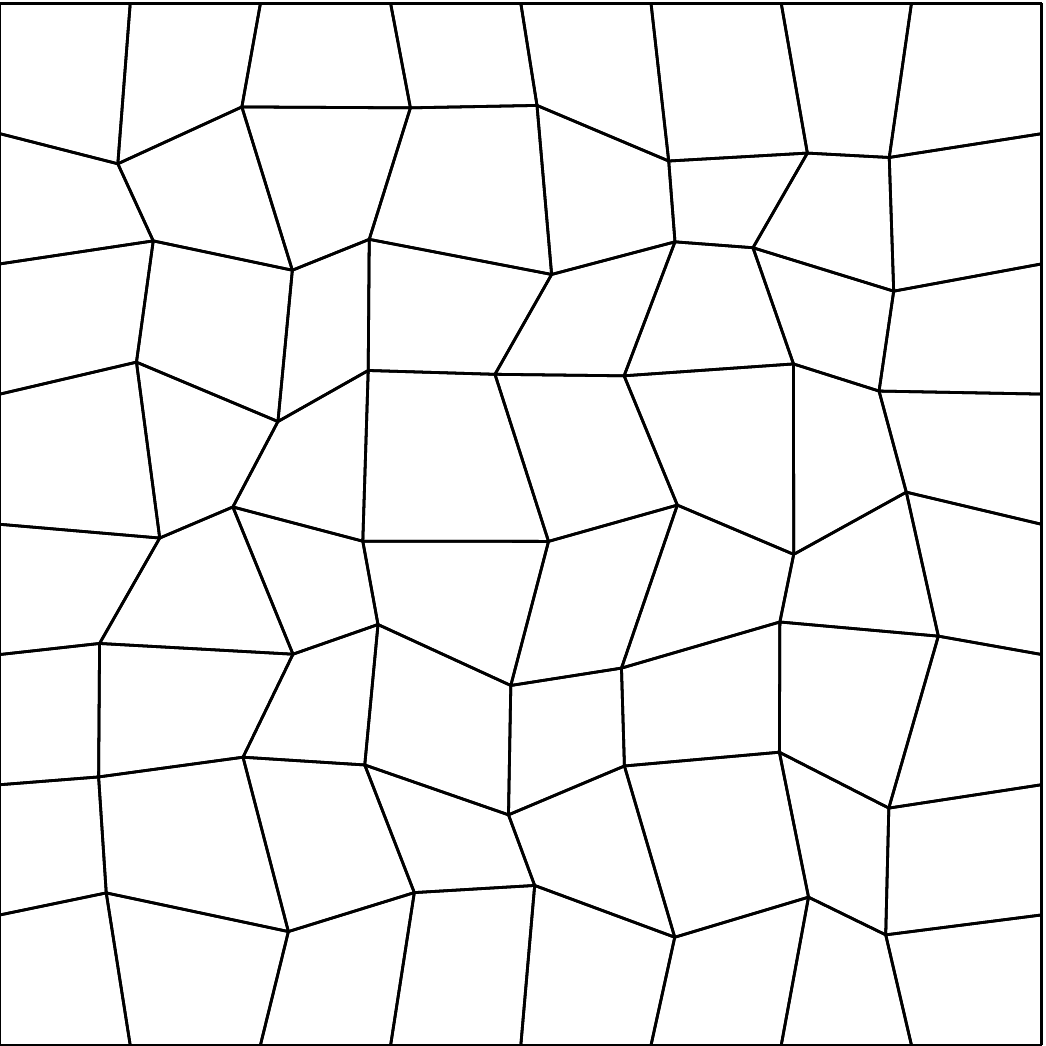}
  \includegraphics[width=0.48\textwidth]{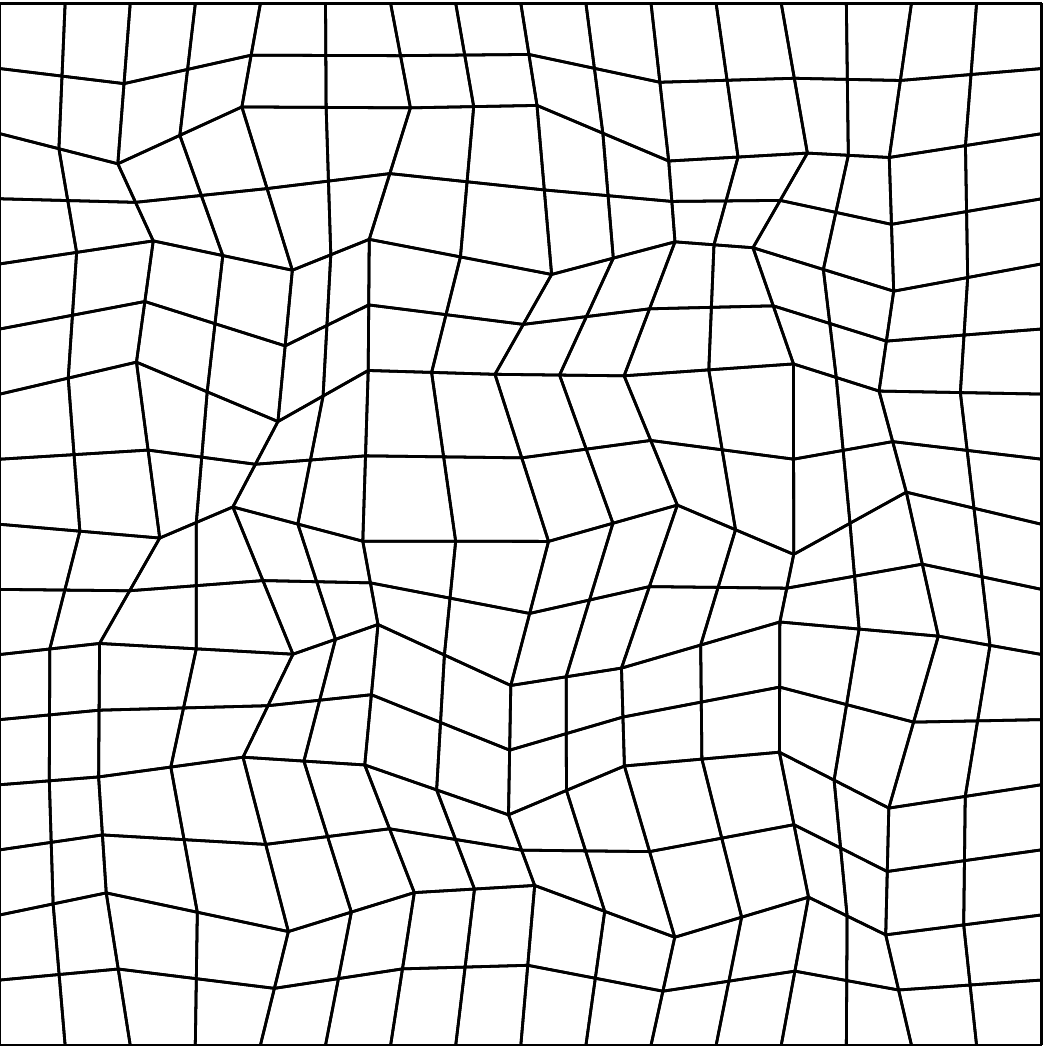}
  \caption{$25\%$ distorted.}
  \label{sfig:gridsdistorted}
  \end{subfigure}
  \begin{subfigure}{0.48\textwidth}
  \includegraphics[width=0.48\textwidth]{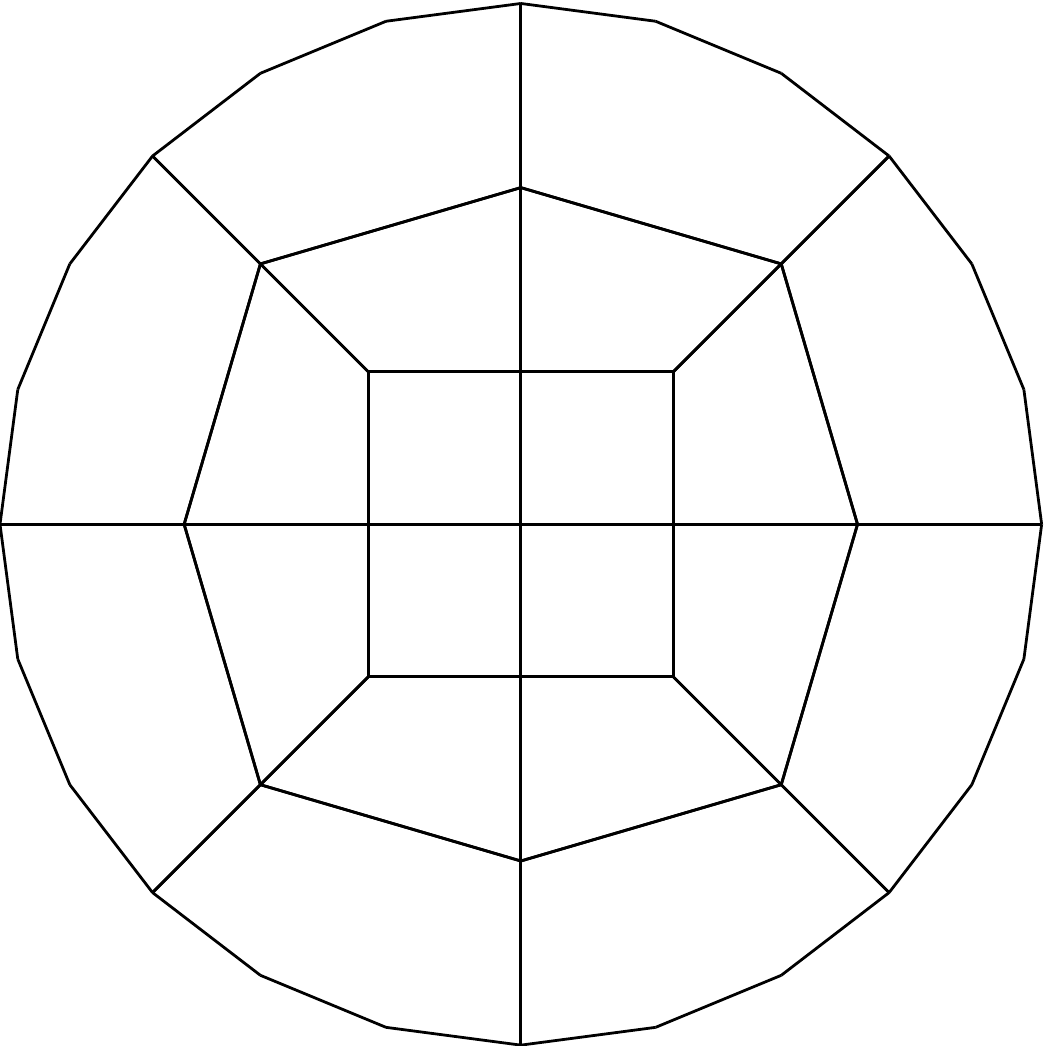}
  \includegraphics[width=0.48\textwidth]{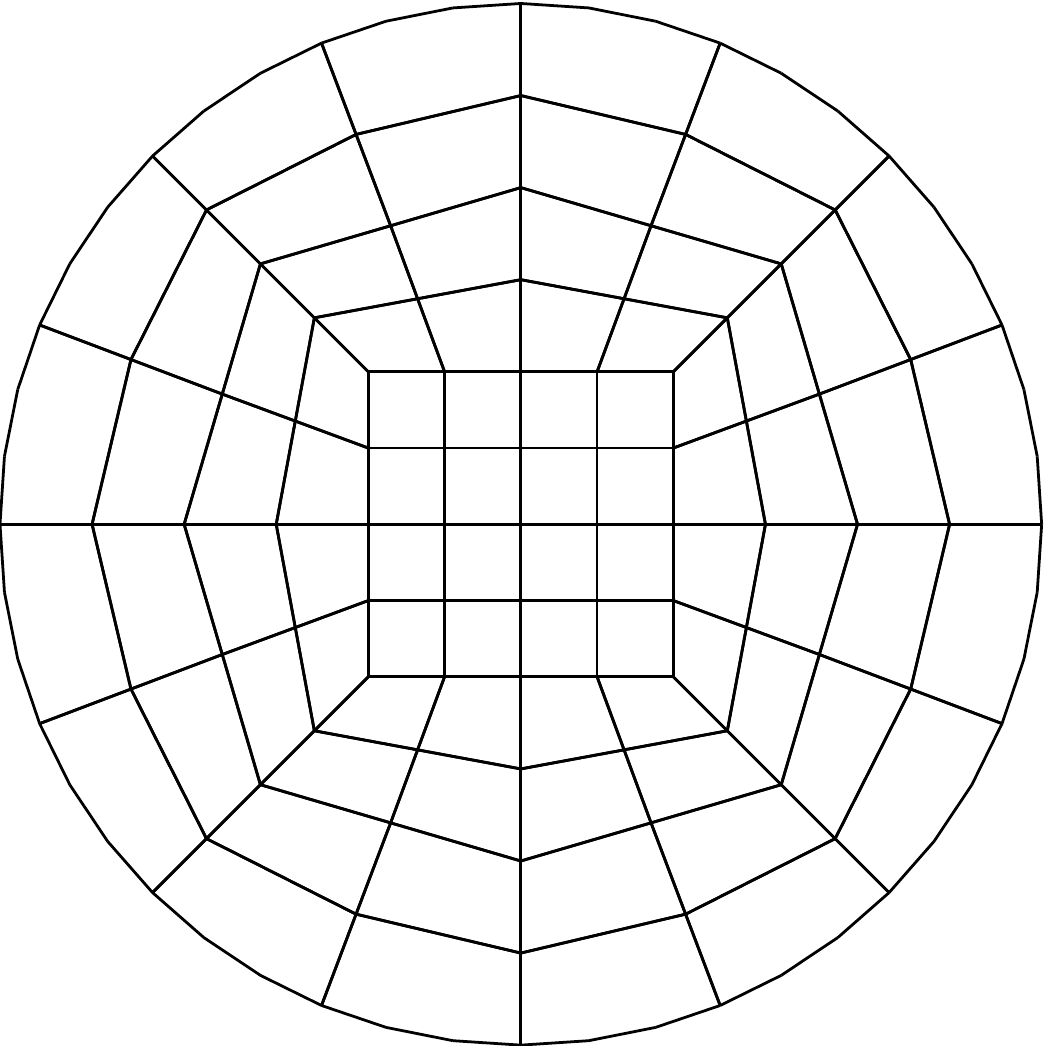}
  \caption{Circular domain.}
  \label{fig:circle}
  \end{subfigure}
  \caption{Initial coarse (left) and subsequent refinement (right) of the distorted and circular grid, respectively.}
  \label{fig:grids}
\end{figure}

We evaluate this method on the meshes in \Cref{fig:grids}, namely a square or a cube with distorted cells and on subdivisions of a circle. For the distorted mesh, we start from a uniform coarse grid of the unit square (cube) with $32 \times 32$ ($8 \times 8 \times 8$) cells. Then, we shift each interior vertex $v$ with respect to the scaling $0.25 h_v$ into a randomly chosen direction, where $h_v$ is the minimal characteristic length of the edges attached to the vertex $v$.
In order to maintain positive definiteness of the interior penalty bilinear form on non-Cartesian meshes, we increase the penalty factor in~\eqref{eq:penalty} to $\ippre = 4$.
In all experiments, we provide fractional iteration counts of the preconditioned conjugate gradient method with a relative accuracy of $10^{-8}$.

\begin{table}[tbh]
\centering
\caption{CG iterations for varying damping $\omega$ on distorted 2D-mesh. Cell based smoothers with one pre- and post-smoothing step for bicubic shape functions.}
\label{tab:cubedist25varydamp}
\begin{tabular}{c cccccc cccccc}
\toprule 
Level L & \multicolumn{6}{c}{ACS} \\
\cmidrule(lr){1-1} \cmidrule(lr){2-7} 
$\omega$ & 0.65 & 0.60 & 0.55 & 0.50 & 0.45 & 0.40 \\
\midrule
3 & \textgreater100 & 54.4 & 40.7 & 38.7 & 42.9 & 48.6  \\
\myrowcolor
4 & \textgreater100 & \textgreater100 & 41.8 & 37.6 & 42.8 & 50.5 \\
5 & \textgreater100 & \textgreater100 & 56.8 & 37.6 & 42.9 & 51.6  \\
\midrule 
Level L & \multicolumn{6}{c}{MCS} \\
\cmidrule(lr){1-1} \cmidrule(lr){2-7} 
$\omega$ & 0.85 & 0.80 & 0.75 & 0.70 & 0.65 & 0.60 \\
\midrule
3 & 23.4 & 23.6 & 24.7 & 25.8 & 27.6 & 29.6 \\
\myrowcolor
4 & \textgreater100 & 28.9 & 23.5 & 24.5 & 25.8 & 27.4 \\
5 & \textgreater100 & 33.9 & 23.3 & 24.3 & 25.5 & 26.6 \\
\bottomrule
\end{tabular}
\end{table}

As predicted by the theory of subspace correction methods, we have to reduce the relaxation parameter $\relaxpre$ in~\eqref{eq:coloredschwarz} compared to Cartesian meshes, see~\Cref{tab:cubedist25varydamp}. There, the best damping factors are $0.50$ for the additive cell based smoother (ACS) and $0.75$ for the multiplicative cell based smoother (MCS). In the following experiments only fractional iteration counts for the best determined damping factor $\relaxpre$ are presented.

\begin{table}[tbh]
\centering
\caption{CG iterations on distorted mesh. Exact (sACS$_{m}$) vs. inexact (ACS$_{m}$) additive cell based smoothers with $m$ pre- and post-smoothing steps.}
\label{tab:cubedist25acp}
\begin{tabular}{c cc cc cc cc}
\toprule 
Level L & \multicolumn{1}{c}{sACS$_{1}$} & \multicolumn{3}{c}{ACS$_{1}$} & \multicolumn{1}{c}{sACS$_{2}$} & \multicolumn{3}{c}{ACS$_{2}$} \\
\cmidrule(lr){1-1} \cmidrule(lr){2-2} \cmidrule(lr){3-5} \cmidrule(lr){6-6} \cmidrule(lr){7-9} 
2D & $k = 3$ & $k = 3$ & $k = 7$ & $k = 15$ & $k = 3$ & $k = 3$ & $k = 7$ & $k = 15$ \\
\cmidrule(lr){1-1} \cmidrule(lr){2-2} \cmidrule(lr){3-5} \cmidrule(lr){6-6} \cmidrule(lr){7-9} 
$\relaxpre$ & 0.70 & 0.50 & 0.50 & 0.50 & 0.70 & 0.50  & 0.50 & 0.50 \\
\midrule
3  & 27.9  & 38.7 & 47.7 & 65.6 & 19.0 & 24.3 & 28.5 & 35.8 \\
\myrowcolor
4  & 27.0  & 37.6 & 47.5 & 65.4 & 18.1 & 23.8 & 28.9 & 36.6 \\
5  & 26.8  & 37.6 & 48.8 & 66.1 & 18.1 & 23.9 & 29.4 & 34.5 \\
\midrule 
3D & $k = 3$ & $k = 3$ & $k = 7$ & $k = 15$ & $k = 3$ & $k = 3$ & $k = 7$ & $k = 15$ \\
\cmidrule(lr){1-1} \cmidrule(lr){2-2} \cmidrule(lr){3-5} \cmidrule(lr){6-6} \cmidrule(lr){7-9} 
$\relaxpre$ & 0.70 & 0.55 & 0.50 & 0.45 & 0.70 & 0.55  & 0.50 & 0.45 \\
\midrule
 1  & 30.0  & 34.4 & 49.7 & 62.6 & 20.4 & 24.5 & 35.9 & 43.3 \\
 \myrowcolor
 2  & 32.4  & 40.3 & 53.7 & 75.6 & 21.1 & 25.9 & 35.9 & 43.3 \\
 3  & 32.0  & 39.5 & 55.7 & 80.5 & 21.3 & 26.1 & 34.8 & 41.3 \\
\bottomrule
\end{tabular}
\end{table}

Starting with the additive cell based smoother, \Cref{tab:cubedist25acp} illustrates that the number of iteration steps increases by less than $\nicefrac12$ when changing from exact local solvers to approximate ones. For two pre- and post-smoothing steps, this increment is even below $\nicefrac14$ for moderate polynomial degrees.
The behavior of this smoother is rather independent of the mesh size $h_L$, but the dependence on polynomial degree seems stronger than on Cartesian meshes in \Cref{tab:cubecp}.
With a second pre- and post-smoothing step each, this dependency becomes less pronounced.
While the gain is only about $\nicefrac23$ for polynomial degree 3, it is almost $\nicefrac12$ for polynomial degree 15.
Therefore, two iteration steps should always be preferred for high polynomial degrees. For lower degrees they are more efficient 
if the transfer operations and/or the coarse grid solver considerably contribute to the computational effort.
The computational cost of computing the average arc lengths~\eqref{eq:avglength} is negligible compared to the setup cost of local solvers.
Hence, the number of arithmetic operations $n_{\textit{FLOP}}$ of setting up the exact and inexact local solvers on non-Cartesian meshes are similar to Cartesian meshes in \Cref{tab:flop-additive} and the smaller number of iterations of $S_\ell^{\textit{exact}}$ (resp. sACS in \Cref{tab:flop-additive}) is outweighed by far by the higher setup cost due to missing tensor structure. 

\begin{table}[tbh]
\centering
\caption{CG iterations on distorted mesh. Inexact multiplicative cell based smoother MCS$_m$ with $m$ pre- and post-smoothing steps and red-black coloring.}
\label{tab:cubedist25mcp}
\begin{tabular}{c cc cc cc}
\toprule 
Level L & \multicolumn{3}{c}{MCS$_{2}$} & \multicolumn{3}{c}{MCS$_{2}$}\\
\cmidrule(lr){1-1} \cmidrule(lr){2-4} \cmidrule(lr){5-7} 
2D & $k = 3$ & $k = 7$ & $k = 15$ & $k = 3$ & $k = 7$ & $k = 15$ \\
\cmidrule(lr){1-1} \cmidrule(lr){2-4} \cmidrule(lr){5-7} 
$\relaxpre$ & 0.75 & 0.50 & 0.50 & 0.75 & 0.50 & 0.50  \\
\midrule
3    & 24.7 & 43.5 & 50.9 & 15.3 & 25.6 & 29.7 \\
\myrowcolor
4    & 23.5 & 42.9 & 49.5 & 14.7 & 24.9 & 28.8 \\
5    & 23.3 & 42.5 & 48.5 & 14.7 & 24.8 & 28.6 \\
\midrule 
3D & $k = 3$ & $k = 7$ & $k = 15$ & $k = 3$ & $k = 7$ & $k = 15$ \\
\cmidrule(lr){1-1} \cmidrule(lr){2-4} \cmidrule(lr){5-7} 
$\relaxpre$ & 0.95 & 0.70 & 0.55 & 0.95 & 0.65 & 0.55 \\
\midrule
1    & 24.7 & 37.7 & 51.7 & 14.4 & 24.8 & 34.5 \\
\myrowcolor
2    & 25.5 & 37.7 & 55.9 & 14.7 & 24.2 & 33.5 \\
3    & 25.8 & 36.5 & 55.8 & 14.8 & 23.9 & 32.7 \\
\bottomrule
\end{tabular}
\end{table} 

Corresponding results for the multiplicative smoother are in \Cref{tab:cubedist25mcp}. Different from \Cref{sec:schwarz}, we use the conjugate gradient method and symmetrize the V-cycle by reverting the order of cells during post-smoothing.
First, we observe that the relaxation parameters are below one, which is again due to the inexact local solvers. Even more, iteration counts are not reduced anymore by one half compared to the additive method, as in \Cref{tab:cubecp}, but only by one third. Thus, due to inexact local solvers, the multiplicative method loses a lot of its attraction.
Similar to the additive smoothers, with two pre- and post-smoothing steps each \Cref{tab:cubedist25mcp} shows an improvement about $\nicefrac{2}{3}$ for tricubic shape functions, about nearly the half for polynomial degree $15$. 

When we treat vertex patches in an analogous way, we observe a much higher sensitivity to mesh distortions. In fact, we were not able to determine reasonable relaxation parameters yielding robust convergence for distortions exceeding 10\%. We attribute this to the fact that the distortions of vertex patches have more degrees of freedom and are thus harder to handle. For the time being, we do not recommend Cartesian surrogates for vertex patches and do not present results.

\begin{table}[tbh]
\centering
\caption{CG iterations on circular mesh. Exact (sACS$_{m}$) vs. inexact (ACS$_{m}$) additive cell based smoothers with $m$ pre- and post-smoothing steps. Entries ``---'' not computed.}
\label{tab:circle-acp}
\begin{tabular}{c cc cc cc cc}
\toprule 
Level L & \multicolumn{1}{c}{sACS$_{1}$} & \multicolumn{3}{c}{ACS$_{1}$} & \multicolumn{1}{c}{sACS$_{2}$} & \multicolumn{3}{c}{ACS$_{2}$} \\
\cmidrule(lr){1-1} \cmidrule(lr){2-2} \cmidrule(lr){3-5} \cmidrule(lr){6-6} \cmidrule(lr){7-9}
2D & $k = 3$ & $k = 3$ & $k = 7$ & $k = 15$ & $k = 3$ & $k = 3$ & $k = 7$ & $k = 15$ \\
\cmidrule(lr){1-1} \cmidrule(lr){2-2} \cmidrule(lr){3-5} \cmidrule(lr){6-6} \cmidrule(lr){7-9} 
$\relaxpre$ & 0.70 & 0.55 & 0.55 & 0.55 & 0.70 & 0.60 & 0.60 & 0.55 \\
\midrule
3  & ---   & ---  & ---  & 57.8 & ---   & ---  & ---  & 33.9 \\
\myrowcolor
4  & ---   & ---  & 44.6 & 58.2 & ---   & ---  & 27.6 & 33.9 \\
5  & 30.5  & 37.2 & 45.6 & 58.8 & 18.4  & 21.6 & 27.5 & 34.7 \\
\myrowcolor
6  & 30.7  & 37.4 & 46.5 & 59.7 & 18.7  & 22.2 & 27.3 & 31.9 \\
7  & 30.9  & 36.6 & 46.8 & ---  & 18.9  & 22.7 & 27.3 & ---  \\
\myrowcolor
8  & 31.5  & 36.3 & ---  & ---  & 19.1  & 22.9 & ---  & ---  \\
\bottomrule
\end{tabular}
\end{table} 

\begin{table}[tbh]
\centering
\caption{CG iterations  on circular mesh. Inexact (MCS$_{m}$) multiplicative cell based smoother with $m$ pre- and post-smoothing steps and coloring by DSATUR algorithm. Entries ``---'' not computed.}
\label{tab:circle-mcp}
\begin{tabular}{c cc cc cc cc}
\toprule 
Level L & \multicolumn{3}{c}{MCS$_{1}$} & \multicolumn{3}{c}{MCS$_{2}$} & Colors\\
\cmidrule(lr){1-1} \cmidrule(lr){2-4} \cmidrule(lr){5-7} \cmidrule(lr){8-8} 
2D & $k = 3$ & $k = 7$ & $k = 15$ & $k = 3$ & $k = 7$ & $k = 15$ \\
\cmidrule(lr){1-1} \cmidrule(lr){2-4} \cmidrule(lr){5-7} 
$\relaxpre$ & 0.95 & 0.70 & 0.65 & 0.95 & 0.70 & 0.65 \\
\midrule
3  & ---   & ---  & 43.9 & ---  & ---  & 25.8 & 4 \\
\myrowcolor
4  & ---   & 33.7 & 45.0 & ---  & 20.6 & 26.5 & 4 \\
5  & 24.7  & 33.7 & 45.7 & 14.6 & 20.8 & 26.9 & 4 \\
\myrowcolor
6  & 25.8  & 32.9 & 45.7 & 14.6 & 20.7 & 25.7 & 4 \\
7  & 25.8  & 33.9 & ---  & 14.6 & 20.4 & ---  & 4 \\
\myrowcolor
8  & 25.8  & ---  & ---  & 14.7 & ---  & ---  & 4 \\
\bottomrule
\end{tabular}
\end{table}

Finally, we performed the same experiments on the subdivisions of the circular domain in \Cref{fig:circle}.  Numerical results for discretizations in two dimensions of mesh level $L$ with $10^5$ to $10^8$ degrees of freedom are shown. The results in \Cref{tab:circle-acp,tab:circle-mcp} exhibit the same characteristics as those for distorted meshes.

\section{Conclusions and outlook}
\label{sec:conclusions}
We have shown in this article that Schwarz smoothers can yield very efficient multigrid solvers if they are implemented by low rank tensor representations. This way, the computational effort for nonoverlapping, cell based smoothing is comparable to that of a matrix-free operator application.
We obtain methods which perform equally well for higher order polynomial spaces. The method, while inherently linked to separable operators and orthogonal mesh cells, can be used to obtain efficient approximate local solvers for more general cases.
By replacing the one-dimensional eigenvalue decompositions of the local solvers by singular value decompositions, it is also applicable to nonsymmetric problems.

The multiplicative smoother with overlapping vertex patches not only outperforms the additive one in terms of iteration steps by far, it is also robust for incompressible vector fields, see for instance~\cite{ArnoldFalkWinther97Hdiv,KanschatMao15}. Albeit it suffers from higher computational cost due to operator applications when switching from one color to the next, it requires less operations than the additive method. There is a computationally more efficient implementation of this smoother, which will reduce its cost even more, in particular for higher order polynomials.

In general, it is not true that low FLOP counts imply fast execution. Nevertheless, we implemented algorithms with similar structure in~\cite{Kronbichler2019Hermite} and showed, that their throughput in MDoF/sec is higher than the one of point smoothers. These in turn are not robust with respect to polynomial degree. In a forthcoming publication, we will focus on performance of a parallel version of the presented work, thus obtaining a fast implementation of a numerically efficient smoother.

\section*{Acknowledgments}
The authors were supported by the German Research Foundation (DFG) under the project 
``High-order discontinuous Galerkin for the exa-scale'' (\mbox{ExaDG}) within the priority program ``Software for Exascale Computing'' (SPPEXA)
and acknowledge support by the state of Baden-W\"urttemberg through bwHPC. The implementation is based on the deal.II library~\cite{dealII91}.

\bibliographystyle{abbrv}
\bibliography{bibexport}
\end{document}